\documentclass{amsart}

\usepackage{amssymb}
\usepackage{amscd}
\usepackage{verbatim}
\usepackage{epsfig}
\begin{document}
\newcommand\Mand{\ \text{and}\ }
\newcommand\Mwith{\ \text{with}\ }
\newcommand\Mfor{\ \text{for}\ }
\newcommand\Mst{\ \text{such that}\ }
\newcommand\Mor{\ \text{or}\ }
\newcommand\Mif{\ \text{if}\ }
\newcommand\Miff{\ \text{iff}\ }
\newcommand\Mthen{\ \text{then}\ }
\newcommand\nin{\notin}
\newcommand\identity{\operatorname{id}}
\newcommand\Id{\operatorname{Id}}
\newcommand\Real{\mathbb{R}}
\newcommand\pos{\Real^+}
\newcommand\Rnp{\Real\setminus\{0\}}
\newcommand\nzero{\setminus\{0\}}
\newcommand\Cx{\mathbb{C}}
\newcommand\Cxp{\Cx^+}
\newcommand\Cxm{\Cx^-}
\newcommand\Nat{\mathbb{N}}
\newcommand\halfNat{{\frac{1}{2}}\mathbb{N}}
\newcommand\intgr{\mathbb{Z}}
\newcommand\im{\operatorname{Im}}
\newcommand\re{\operatorname{Re}}
\newcommand\sign{\operatorname{sign}}
\newcommand\codim{\operatorname{codim}}
\newcommand\End{\operatorname{End}}
\newcommand\Ker{\operatorname{Ker}}
\newcommand\Hom{\operatorname{Hom}}
\newcommand\tr{\operatorname{tr}}
\newcommand\Tr{\operatorname{Tr}}
\newcommand\ideal{{\mathcal I}}
\newcommand\Span{\operatorname{span}}
\newcommand\image{\operatorname{image}}
\newcommand\Range{\operatorname{Ran}}
\newcommand\Graph{\operatorname{graph}}
\newcommand\slim{\operatornamewithlimits{s-lim}}
\newcommand\spp{\operatorname{sp}}
\newcommand\sll{\operatorname{sl}}
\newcommand\sol{\operatorname{so}}
\newcommand\SL{\operatorname{SL}}
\newcommand\SO{\operatorname{SO}}
\newcommand\On{\operatorname{O}}
\newcommand\pa{\partial}
\newcommand\Rn{\Real^n}
\newcommand\Rm{\Real^m}
\newcommand\RN{\Real^N}
\newcommand\RtN{\Real^{2N}}
\newcommand\RM{\Real^M}
\newcommand\sphere{\mathbb{S}}
\newcommand\Sn{\sphere^{n-1}}
\newcommand\Sm{\sphere^{m-1}}
\newcommand\Snp{\sphere^n_+}
\newcommand\Smp{\sphere^m_+}
\newcommand\SN{\sphere^{N-1}}
\newcommand\SNp{\sphere^N_+}
\newcommand\circlep{\sphere^1_+}
\newcommand\Phom{P_{h}}
\newcommand\Shom{S_{h}}
\newcommand\distance{\operatorname{dist}}
\newcommand\cl{\operatorname{cl}}
\newcommand\interior{\operatorname{int}}
\newcommand\Fa{\operatorname{Fa}}
\newcommand\ff{\operatorname{ff}}
\newcommand\mf{\operatorname{mf}}
\newcommand\cf{\operatorname{cf}}
\newcommand\scf{\operatorname{sf}}
\newcommand\lf{\operatorname{lf}}
\newcommand\rf{\operatorname{rf}}
\newcommand\indfam{{\mathcal K}}
\newcommand\fraka{{\mathfrak a}}
\newcommand\calE{{\mathcal E}}
\newcommand\calA{{\mathcal A}}
\newcommand\calB{{\mathcal B}}
\newcommand\calH{{\mathcal H}}
\newcommand\calR{{\mathcal R}}
\newcommand\calO{{\mathcal O}}
\newcommand\calJ{{\mathcal J}}
\newcommand\calM{{\mathcal M}}
\newcommand\calN{{\mathcal N}}
\newcommand\calX{{\mathcal X}}
\newcommand\calF{{\mathcal F}}
\newcommand\calG{{\mathcal G}}
\newcommand\calT{{\mathcal T}}
\newcommand\calP{{\mathcal P}}
\newcommand\calPt{\widetilde{\mathcal P}}
\newcommand\calS{{\mathcal S}}
\newcommand\calSt{\widetilde{\mathcal S}}
\newcommand\calC{{\mathcal C}}
\newcommand\calCt{{\tilde {\mathcal C}}}
\newcommand\calCL{{\mathcal C}_{\text L}}
\newcommand\calCR{{\mathcal C}_{\text R}}
\newcommand\Cinf{{\mathcal C}^{\infty}}
\newcommand\dist{{\mathcal C}^{-\infty}}
\newcommand\dCinf{\dot\Cinf}
\newcommand\ddist{\dot\dist}
\newcommand\Cj{{\mathcal C}^j}
\newcommand\Linf{L^{\infty}}
\newcommand\phg{{\text{phg}}}
\newcommand\comp{{\text{comp}}}
\newcommand\bcon{{\mathcal A}}
\newcommand\bconc{{\mathcal A}_{\text{phg}}}
\newcommand\Sch{{\mathcal S}}
\newcommand\temp{\Sch^{\prime}}
\newcommand\Diff{\operatorname{Diff}}
\newcommand\Diffb{\operatorname{Diff}_{\text{b}}}
\newcommand\Diffc{\operatorname{Diff}_{\text{c}}}
\newcommand\Diffsc{\operatorname{Diff}_{\text{sc}}}
\newcommand\DiffI{\operatorname{Diff}_{\text{I}}}
\newcommand\DiffIq{\operatorname{Diff}_{\text{I},q}}
\newcommand\sing{\text{sing}}
\newcommand\reg{\text{reg}}
\newcommand\supp{\operatorname{supp}}
\newcommand\ssupp{\operatorname{sing\ supp}}
\newcommand\csupp{\operatorname{cone\ supp}}
\newcommand\esupp{\operatorname{ess\ supp}}
\newcommand\Fr{{\mathcal F}}
\newcommand\Frinv{\Fr^{-1}}
\newcommand\bop{{\mathcal B}}
\newcommand\spec{\operatorname{spec}}
\newcommand\pspec{\spec_{pp}}
\newcommand\cspec{\spec_{c}}
\newcommand\FIO{{\mathcal I}}
\newcommand\SP{\operatorname{RC}}
\newcommand\RC{\operatorname{RC}}
\newcommand\Symc{S_c}
\newcommand\Symca{S_c^{\alpha}}
\newcommand\Symczero{S_c^{0,...,0}}
\newcommand\sci{{}^{\text{sc}}}
\newcommand\sct{\sci T^*}
\newcommand\scT{\sci T}
\newcommand\scdt{\sci \dot T^*}
\newcommand\dS{\dot S^*}
\newcommand\dT{\dot T^*}
\newcommand\dSreg{\dot\Sigma_{\text reg}}
\newcommand\scct{\sci\bar{T}^*}
\newcommand\Csc{C_{\text{sc}}}
\newcommand\SNpscd{(\SNp)^2_{\text{sc}}}
\newcommand\scdiag{\Delta_{\text{sc}}}
\newcommand\projscl{\pi^L_{\text{sc}}}
\newcommand\projscr{\pi^R_{\text{sc}}}
\newcommand\scHL{\sci H^{2,0}_{|\zeta|^2-\lambda^2}}
\newcommand\scHrg{\sci H^{2,0}_{\sqrt{g}}}
\newcommand\Hsc{H_{\text{sc}}}
\newcommand\WF{\operatorname{WF}}
\newcommand\WFp{\operatorname{WF^{\prime}}}
\newcommand\WFsc{\operatorname{WF}_{\text{sc}}}
\newcommand\WFscp{\operatorname{WF_{sc}^{\prime}}}
\newcommand\WFC{\operatorname{WF}_C}
\newcommand\WFCi{\operatorname{WF}_{C_i}}
\newcommand\elliptic{\operatorname{ell}}
\newcommand\Psop{\operatorname{\Psi}}
\newcommand\Psiscrs{\operatorname{\Psi_{sc}^{-2,\infty}}}
\newcommand\Psiscr{\operatorname{\Psi_{sc}^{-2,0}}}
\newcommand\Psiscrm{\operatorname{\Psi_{sc}^{0,2}}}
\newcommand\PsiscHam{\operatorname{\Psi_{sc}^{2,0}}}
\newcommand\Psisci{\operatorname{\Psi_{sc}^{*,*}}}
\newcommand\Psiscid{\operatorname{\Psi_{sc}^{0,0}}}
\newcommand\Psiscis{\operatorname{\Psi_{sc}^{0,\infty}}}
\newcommand\Psiscsi{\operatorname{\Psi_{sc}^{-\infty,0}}}
\newcommand\Psiscs{\operatorname{\Psi_{sc}^{-\infty,\infty}}}
\newcommand\Psiscalg{\operatorname{\Psi_{sc}^{\infty,-\infty}}}
\newcommand\nullHam{{\mathcal N}}
\newcommand\charD{\Sigma_{\Delta-\lambda^2}}
\newcommand\charLap{\Sigma_{\Delta-\lambda}}
\newcommand\Snl{\Sn_{\lambda}}
\newcommand\SNl{\SN_{\lambda}}
\newcommand\gammat{\tilde\gamma}
\newcommand\gammasc{\gamma}
\newcommand\Tau{\mathcal{T}}
\newcommand\taut{\tilde\tau}
\newcommand\taub{\bar\tau}
\newcommand\Nout{N^+_{\lambda}}
\newcommand\Nin{N^-_{\lambda}}
\newcommand\Nio{N^{\pm}_{\lambda}}
\newcommand\El{E_{\lambda}}
\newcommand\Elt{\tilde E_{\lambda}}
\newcommand\Eil{E^i_{\lambda}}
\newcommand\Ejl{E^j_{\lambda}}
\newcommand\Eajl{E^{\alpha_j}_{\lambda}}
\newcommand\Eilt{\tilde E^i_{\lambda}}
\newcommand\Np{N^+}
\newcommand\Nm{N^-}
\newcommand\Npm{N^{\pm}}
\newcommand\Fin{F^-(\lambda)}
\newcommand\Fini{F^-_i(\lambda)}
\newcommand\Fout{F^+(\lambda)}
\newcommand\Fouti{F^+_i(\lambda)}
\newcommand\Foutj{F^+_j(\lambda)}
\newcommand\Rout{R^+_{\lambda}}
\newcommand\Routl{R^+_{\lambda^2}}
\newcommand\Routsgnl{R^{\sign\lambda}_{\lambda^2}}
\newcommand\Rin{R^-_{\lambda}}
\newcommand\Rinl{R^-_{\lambda^2}}
\newcommand\Rinsgnl{R^{-\sign\lambda}_{\lambda^2}}
\newcommand\Rio{R^{\pm}_{\lambda}}
\newcommand\Riol{R^{\pm}_{\lambda^2}}
\newcommand\Roi{R^{\mp}_{\lambda}}
\newcommand\Roil{R^{\mp}_{\lambda^2}}
\newcommand\Riob{R^{\pm}}
\newcommand\Roib{R^{\mp}}
\newcommand\Tio{T^{\pm}}
\newcommand\Tiob{T^{\pm}_{\ff}}
\newcommand\Toi{T^{\mp}}
\newcommand\Toib{T^{\mp}_{\ff}}
\newcommand\TIiob{T_I^{\pm}}
\newcommand\Rinb{R^-}
\newcommand\Rinbsgnl{R^{-\sign\lambda}}
\newcommand\Tin{T^-}
\newcommand\Tinb{T^-_{\ff}}
\newcommand\TIinb{T^-_I}
\newcommand\Routb{R^+}
\newcommand\Routbsgnl{R^{\sign\lambda}}
\newcommand\Tout{T^+}
\newcommand\Toutb{T^+_{\ff}}
\newcommand\TIoutb{T^+_I}
\newcommand\Rlkf{(|\xib|^2-(\lambda-i0)^2)^{-1}}
\newcommand\Rlk{\rho_0(\lambda)}
\newcommand\Rmlk{\rho_0(-\lambda)}
\newcommand\Rpmlk{\rho_0(\pm\lambda)}
\newcommand\Rlka{\rho_1(\lambda)}
\newcommand\Rlkb{\rho_2(\lambda)}
\newcommand\Rilk{\rho_i(\lambda)}
\newcommand\reduced{\natural}
\newcommand\Rlf{R_0(\lambda)}
\newcommand\Rla{R_1(\lambda)}
\newcommand\Rlb{R_2(\lambda)}
\newcommand\Ril{R_i(\lambda)}
\newcommand\Rlj{R_j(\lambda)}
\newcommand\Rlft{R_0(\lambda)}
\newcommand\Rflambda{R_0^{\reduced}(\sigma)}
\newcommand\RV{R^{\reduced}_V}
\newcommand\Rfsigma{R_0^{\reduced}(\sigma)}
\newcommand\Rfsigmah{R_0^{\reduced}(\sigma^{1/2})}
\newcommand\Rfzero{R_0^{\reduced}(0)}
\newcommand\RlV{R^{\reduced}_V(\sigma)}
\newcommand\RlVi{R^{\reduced}_{V_i}(\sigma)}
\newcommand\RlVt{R_V(\lambda)}
\newcommand\RlVtL{{R}_V^L(\lambda)}
\newcommand\RlVtR{{R}_V^R(\lambda)}
\newcommand\RlVit{{R}_{V_i}(\lambda)}
\newcommand\RlVta{{R}_V^{(1)}(\lambda)}
\newcommand\RlVtk{{R}_V^{(k)}(\lambda)}
\newcommand\RlVatV{{R}_{V_{\alpha}}(\lambda)V_{\alpha}}
\newcommand\RlVatVa{{R}_{V_{\alpha_1}}(\lambda)V_{\alpha_1}}
\newcommand\RlVatVb{{R}_{V_{\alpha_2}}(\lambda)V_{\alpha_2}}
\newcommand\RlVatVk{{R}_{V_{\alpha_k}}(\lambda)V_{\alpha_k}}
\newcommand\RlVatVkk{{R}_{V_{\alpha_{k+1}}}(\lambda)V_{\alpha_{k+1}}}
\newcommand\RlVaptV{{R}_{V_{\alpha'}}(\lambda)V_{\alpha'}}
\newcommand\RlVapptV{{R}_{V_{\alpha''}}(\lambda)V_{\alpha''}}
\newcommand\RlVajtV{{R}_{V_{\alpha_j}}(\lambda)V_{\alpha_j}}
\newcommand\RlVaktV{{R}_{V_{\alpha_k}}(\lambda)V_{\alpha_k}}
\newcommand\RlVakktV{{R}_{V_{\alpha_{k+1}}}(\lambda)V_{\alpha_{k+1}}}
\newcommand\Tl{T(\lambda)}
\newcommand\Tlt{\tilde\Tl}
\newcommand\Tltp{\tilde T'(\lambda)}
\newcommand\Tltpp{\tilde T''(\lambda)}
\newcommand\Tli{T_i(\lambda)}
\newcommand\Tlit{\tilde\Tli}
\newcommand\Tlip{T_i'(\lambda)}
\newcommand\Tlipp{T_i''(\lambda)}
\newcommand\Tlj{T_j(\lambda)}
\newcommand\Tla{T_{\alpha}(\lambda)}
\newcommand\Tlaa{T_{\alpha_1}(\lambda)}
\newcommand\Tlab{T_{\alpha_2}(\lambda)}
\newcommand\Tlak{T_{\alpha_k}(\lambda)}
\newcommand\Tlakt{\tilde\Tlak}
\newcommand\Tlaj{T_{\alpha_j}(\lambda)}
\newcommand\Tlajj{T_{\alpha_{j+1}}(\lambda)}
\newcommand\Tlajp{T_{\alpha_j}'(\lambda)}
\newcommand\Tlajpt{\tilde\Tlajp}
\newcommand\Tlajt{\tilde\Tlaj}
\newcommand\Tlakk{T_{\alpha_{k+1}}(\lambda)}
\newcommand\Tlakkp{T_{\alpha_{k+1}}'(\lambda)}
\newcommand\Tlap{T_{\alpha'}(\lambda)}
\newcommand\Tlapt{\tilde\Tlap}
\newcommand\Tlapp{T_{\alpha''}(\lambda)}
\newcommand\Tkl{T^{(k)}(\lambda)}
\newcommand\Tcl{T^{\flat}(\lambda)}
\newcommand\Fl{F(\lambda)}
\newcommand\BlVt{\tilde B_V(\lambda)}
\newcommand\KBlVt{K_{\BlVt}}
\newcommand\BlVaat{B_{V_{\alpha_1}}(\lambda)}
\newcommand\BV{B_V}
\newcommand\Bone{B_1}
\newcommand\Btwo{B_2}
\newcommand\Bthree{B_3}
\newcommand\Banyj{B_j}
\newcommand\PlV{P_V(\lambda)}
\newcommand\PlVc{P_V^{\flat}(\lambda)}
\newcommand\Pl{P_0(\lambda)}
\newcommand\SVl{S_V(\lambda)}
\newcommand\Sjr{S_j^{\reduced}}
\newcommand\Rkp{{\mathcal R}^k_+}
\newcommand\Rkm{{\mathcal R}^k_-}
\newcommand\Rkpm{{\mathcal R}^k_{\pm}}
\newcommand\Phys{{\mathcal P}}
\newcommand\Pc{\overline{\mathcal P}}
\newcommand\pip{\pi^{\perp}}
\newcommand\pipa{\pi_\partial}
\newcommand\gammapa{\gamma_\partial}
\newcommand\pipah{\hat\pi_\partial}
\newcommand\pit{\tilde\pi}
\newcommand\xit{\tilde\xi}
\newcommand\zetat{\tilde\zeta}
\newcommand\etat{\tilde\eta}
\newcommand\sigmat{\tilde\sigma}
\newcommand\sigmahat{\hat\sigma}
\newcommand\thetat{\tilde\theta}
\newcommand\psit{\tilde\psi}
\newcommand\phit{\tilde\phi}
\newcommand\chit{\tilde\chi}
\newcommand\rhot{\tilde\rho}
\newcommand\xib{\bar\xi}
\newcommand\zetab{\bar\zeta}
\newcommand\thetab{\bar\theta}
\newcommand\etab{\bar\eta}
\newcommand\iotal{\iota_{\lambda}}
\newcommand\rhoat{\rhot_{\alpha_1}}
\newcommand\Lambdat{\tilde\Lambda}
\newcommand\Lambdati{\tilde\Lambda^{\text{in}}}
\newcommand\Lambdato{\tilde\Lambda^{\text{out}}}
\newcommand\Lambdatp{\tilde\Lambda^{\text{prop}}}
\newcommand\Lambdai{\Lambda^{\text{in}}}
\newcommand\Lambdao{\Lambda^{\text{out}}}
\newcommand\poles{\Lambda'}
\newcommand\rpoles{\Lambda_p}
\newcommand\thresholds{\Lambda}
\newcommand\Vt{\tilde V}
\newcommand\It{\tilde I}
\newcommand\half{{\frac{1}{2}}}
\newcommand\sigmah{\sigma^{1/2}}
\newcommand\bX{\partial X}
\newcommand\bXb{\partial \Xb}
\newcommand\Deltabt{\tilde\Delta_0}
\newcommand\strip{\Omega_T}
\newcommand\Kf{K^{\flat}}
\newcommand\Gs{G^{\sharp}}
\newcommand\Gt{\tilde G}
\newcommand\Osc{\sci\Omega}
\newcommand\OSc{{}^\Scl\Omega}
\newcommand\Osch{\sci\Omega^{\half}}
\newcommand\Oscmh{\sci\Omega^{-\half}}
\newcommand\Isc{I_{sc}}
\newcommand\os{{\text{os}}}
\newcommand\Qzl{Q^0_{-\lambda}}
\newcommand\Lie{{\mathcal L}}
\newcommand\bl{{\text b}}
\newcommand\scl{{\text{sc}}}
\newcommand\sccl{{\text{scc}}}
\newcommand\Scl{{\text{Sc}}}
\newcommand\ScLl{{\text{Sc,L}}}
\newcommand\ScRl{{\text{Sc,R}}}
\newcommand\Sccl{{\text{Scc}}}
\newcommand\sus{{\text{sus}}}
\newcommand\ssl{{\text{ss}}}
\newcommand\XXb{X^2_\bl}
\newcommand\XXbt{\Xt^2_\bl}
\newcommand\XXsc{X^2_\scl}
\newcommand\XXsct{\Xt^2_\scl}
\newcommand\XXSc{X^2_\Scl}
\newcommand\XXSct{\Xt^2_\Scl}
\newcommand\XXScL{X^2_\ScLl}
\newcommand\XXScR{X^2_\ScRl}
\newcommand\MMsc{M^2_\scl}
\newcommand\Deltab{\Delta_\bl}
\newcommand\Deltasc{\Delta_\scl}
\newcommand\DeltaSc{\Delta_\Scl}
\newcommand\DeltaScL{\Delta_\ScLl}
\newcommand\DeltaScR{\Delta_\ScRl}
\newcommand\prs{\sigma}
\newcommand\Nsc{N_\scl}
\newcommand\Nscp{N_{\scl,p}}
\newcommand\Nff{N_{\ff}}
\newcommand\Nffz{N_{\ff,0}}
\newcommand\Nffzp{N_{\ff,0,p}}
\newcommand\Nffl{N_{\ff,l}}
\newcommand\Nffml{N_{\ff,-l}}
\newcommand\Nmf{N_{\mf}}
\newcommand\Nmfz{N_{\mf,0}}
\newcommand\Nmfl{N_{\mf,l}}
\newcommand\Nmfml{N_{\mf,-l}}
\newcommand\ffb{\operatorname{bf}}
\newcommand\Ffb{\operatorname{bf'}}
\newcommand\ffsc{\operatorname{sf}}
\newcommand\ffSc{\operatorname{sf_C}}
\newcommand\Ffsc{\operatorname{sf'}}
\newcommand\rff{\rho_{\ff}}
\newcommand\rmf{\rho_{\mf}}
\newcommand\rffb{\rho_{\ffb}}
\newcommand\rffsc{\rho_{\ffsc}}
\newcommand\rFfsc{\rho_{\Ffsc}}
\newcommand\rffSc{\rho_{\ffSc}}
\newcommand\rinf{\rho_{\infty}}
\newcommand\CL{C_L}
\newcommand\CR{C_R}
\newcommand\betab{\beta_\bl}
\newcommand\betasc{\beta_\scl}
\newcommand\betaSc{\beta_\Scl}
\newcommand\BetaSc{\bar\beta_\Scl}
\newcommand\betaScL{\beta_\ScLl}
\newcommand\betaScR{\beta_\ScRl}
\newcommand\ScT{{}^\Scl T^*}
\newcommand\SccT{{}^\Scl \bar T^*}
\newcommand\ScS{{}^\Scl S^*}
\newcommand\Tb{{}^\bl T}
\newcommand\Tsc{{}^\scl T}
\newcommand\TSc{{}^\Scl T}
\newcommand\CSc{C_\Scl}
\newcommand\Lambdasc{{}^\scl\Lambda}
\newcommand\XXXb{X^3_\bl}
\newcommand\XXXsc{X^3_\scl}
\newcommand\XXXSc{X^3_\Scl}
\newcommand\XXXScO{X^3_{\Scl,O}}
\newcommand\XXXScF{X^3_{\Scl,F}}
\newcommand\XXXScS{X^3_{\Scl,S}}
\newcommand\XXXScC{X^3_{\Scl,C}}
\newcommand\KDsc{\operatorname{KD^{\half}_\scl}}
\newcommand\KDSc{\operatorname{KD^{\half}_\Scl}}
\newcommand\KDScEF{\operatorname{KD^{E,F}_\Scl}}
\newcommand\Oh{\operatorname{\Omega^{\half}}}
\newcommand\WFSc{\WF_\Scl}
\newcommand\WFtSc{\WF_{\text 3sc}}
\newcommand\WFScmf{\WF_{\Scl,\mf}}
\newcommand\WFScff{\WF_{\Scl,\ff}}
\newcommand\WFScs{\WF_{\Scl,\prs}}
\newcommand\WFScp{\WF'_\Scl}
\newcommand\WFScmfp{\WF'_{\Scl,\mf}}
\newcommand\WFScffp{\WF'_{\Scl,\ff}}
\newcommand\WFScsp{\WF'_{\Scl,\prs}}
\newcommand\Diffscc{\Diff_\sccl}
\newcommand\DiffSc{\Diff_\Scl}
\newcommand\DiffScc{\Diff_\Sccl}
\newcommand\DiffscI{\Diff_{\scl,\text{I}}}
\newcommand\VscI{\Vf_{\scl,\text{I}}}
\newcommand\DiffsV{\operatorname{Diff}_{\sus(V)}}
\newcommand\DiffsVsc{\operatorname{Diff}_{\sus(V),\scl}}
\newcommand\DiffsVCsc{\operatorname{Diff}_{\sus(V)-C,\scl}}   
\newcommand\Psisc{\Psop_\scl}
\newcommand\Psiscc{\Psop_\sccl}
\newcommand\Psiss{\Psop_\ssl}
\newcommand\Psisch{\Psop_{\scl,h}}
\newcommand\Psiscch{\Psop_{\sccl,h}}
\newcommand\PsiSc{\Psop_\Scl}
\newcommand\PsiScph{\Psop_{\Scl,\phi}}
\newcommand\PsiScra{\Psop_{\Scl,\rho^\sharp_a}}
\newcommand\PsiScc{\Psop_\Sccl}
\newcommand\PsiSccml{\Psop^{m,l}_\Sccl}
\newcommand\PsiScxx{\Psop^{*,*}_\Scl}
\newcommand\PsiScml{\Psop^{m,l}_\Scl}
\newcommand\PsiScmz{\Psop^{m,0}_\Scl}
\newcommand\PsiScmmz{\Psop^{-m,0}_\Scl}
\newcommand\PsiSckz{\Psop^{k,0}_\Scl}
\newcommand\PsiScmmml{\Psop^{-m,-l}_\Scl}
\newcommand\Psiscmkk{\Psop^{-k,k}_\scl}
\newcommand\Psiscmmmkk{\Psop^{-m-k,k}_\scl}
\newcommand\Psiscmoo{\Psop^{-1,1}_\scl}
\newcommand\Psiscmz{\Psop^{m,0}_\scl}
\newcommand\Psiscmmz{\Psop^{-m,0}_\scl}
\newcommand\PsiSckmkl{\Psop^{km,kl}_\Scl}
\newcommand\PsiScmplp{\Psop^{m',l'}_\Scl}
\newcommand\PsiScmmpllp{\Psop^{m+m',l+l'}_\Scl}
\newcommand\Psiscml{\Psop^{m,l}_\scl}
\newcommand\PsiScid{\Psop^{0,0}_\Scl}
\newcommand\PsiSczo{\Psop^{0,1}_\Scl}
\newcommand\PsiScmii{\Psop^{-\infty,\infty}_\Scl}
\newcommand\PsiScmiz{\Psop^{-\infty,0}_\Scl}
\newcommand\PsiScmoo{\Psop^{-1,1}_\Scl}
\newcommand\PsisCid{\Psop^{0,0}_{\scl-C}}
\newcommand\PsisC{\Psop_{\scl-C}}
\newcommand\Psiinf{\Psop_{\infty}}
\newcommand\Psiinfid{\Psop_{\infty}^0}
\newcommand\PsiFinf{\Psop_{\infty-\Fr}}
\newcommand\PsisVscml{\Psop^{m,l}_{\sus(V),\scl}}
\newcommand\PsisVsc{\Psop_{\sus(V),\scl}}
\newcommand\PsisVpsc{\Psop_{\sus(V_p),\scl}}
\newcommand\PsisVCSc{\Psop_{\sus(V)-C,\scl}}
\newcommand\SFinf{S_{\infty-\Fr}}
\newcommand\YsVC{Y^2_{\sus(V)-C,\scl}}
\newcommand\ffYsc{\ffsc_{\sus(V)}}
\newcommand\SXC{S(X;C)}
\newcommand\Ios{I_{\text{os}}}
\newcommand\pbL{\pi^2_{\bl,{\text L}}}
\newcommand\pbR{\pi^2_{\bl,{\text R}}}
\newcommand\pscL{\pi^2_{\scl,{\text L}}}
\newcommand\pscR{\pi^2_{\scl,{\text R}}}
\newcommand\PbO{\pi^3_{\bl,{\text O}}}
\newcommand\PscO{\pi^3_{\scl,{\text O}}}
\newcommand\PScO{\pi^3_{\Scl,{\text O}}}
\newcommand\PScF{\pi^3_{\Scl,{\text F}}}
\newcommand\PScC{\pi^3_{\Scl,{\text C}}}
\newcommand\PScS{\pi^3_{\Scl,{\text S}}}
\newcommand\pScL{\pi^2_{\Scl,{\text L}}}
\newcommand\pScR{\pi^2_{\Scl,{\text R}}}
\newcommand\CLF{\CL^F}
\newcommand\CLO{\CL^O}
\newcommand\CLS{\CL^S}
\newcommand\CLC{\CL^C}
\newcommand\DeltaYb{\Delta_{\bl,Y}}
\newcommand\DeltaYsc{\Delta_{\sus-\scl}}
\newcommand\diag{\operatorname{diag}}
\newcommand\Vf{{\mathcal V}}
\newcommand\Vb{{\mathcal V}_{\bl}}
\newcommand\Vsc{{\mathcal V}_{\scl}}
\newcommand\VSc{{\mathcal V}_{\Scl}}
\newcommand\VfI{\Vf_{\text{I}}}
\newcommand\VfIq{\Vf_{\text{I},q}}
\newcommand\scH{{}^\scl H}
\newcommand\scHg{\scH_g}
\newcommand\Hss{H_\ssl}
\newcommand\xh{\hat x}
\newcommand\yh{\hat y}
\newcommand\sh{\hat s}
\newcommand\rh{\hat r}
\newcommand\Yh{\hat Y}
\newcommand\Zh{\hat Z}
\newcommand\Yb{\bar Y}
\newcommand\hb{\bar h}
\newcommand\xih{\hat\xi}
\newcommand\etah{\hat\eta}
\newcommand\muh{\hat\mu}
\newcommand\mub{\bar\mu}
\newcommand\nub{\bar\nu}
\newcommand\mubh{\widehat{\bar\mu}}
\newcommand\yb{\bar y}
\newcommand\zb{\bar z}
\newcommand\ub{\bar u}
\newcommand\Qb{\bar Q}
\newcommand\Wbp{{\bar W}^\perp}
\newcommand\Wp{W^\perp}
\newcommand\Kt{\tilde K}
\newcommand\Wt{\tilde W}
\newcommand\Ut{\tilde U}
\newcommand\yt{\tilde y}
\newcommand\ut{\tilde u}
\newcommand\vt{\tilde v}
\newcommand\ft{\tilde f}
\newcommand\htil{\tilde h}
\newcommand\St{\tilde S}
\newcommand\Pt{\tilde P}
\newcommand\Rt{\tilde R}
\newcommand\qt{\tilde q}
\newcommand\Qt{\tilde Q}
\newcommand\Xb{\bar X}
\newcommand\lambdat{\tilde\lambda}
\newcommand\betat{\tilde\beta}
\newcommand\Phit{\tilde\Phi}
\newcommand\epst{\tilde\epsilon}
\newcommand\ep{\epsilon}
\newcommand\bt{\tilde b}
\newcommand\Xt{\tilde X}
\newcommand\Mt{\tilde M}
\newcommand\At{\tilde A}
\newcommand\Et{\tilde E}
\newcommand\Ht{\tilde H}
\newcommand\at{\tilde a}
\newcommand\Ct{\tilde C}
\newcommand\pih{\hat\pi}
\newcommand\Rh{\hat R}
\newcommand\Ah{\hat A}
\newcommand\Bh{\hat B}
\newcommand\Ch{\hat C}
\newcommand\Gh{\hat G}
\newcommand\Hh{\hat H}
\newcommand\Qh{\hat Q}
\newcommand\Ph{\hat P}
\newcommand\Nh{\hat N}
\newcommand\Sh{\hat S}
\newcommand\Gcal{{\mathcal G}}
\newcommand\GcalC{{\mathcal G}_C}
\newcommand\Jcal{{\mathcal J}}
\newcommand\JcalC{{\mathcal J}_C}
\newcommand\evpr{\lambda_1}
\newcommand\evth{\lambda_0}
\setcounter{secnumdepth}{3}
\newtheorem{lemma}{Lemma} %[section]
\newtheorem{prop}[lemma]{Proposition}
\newtheorem{thm}[lemma]{Theorem}
\newtheorem{cor}[lemma]{Corollary}
\newtheorem{result}[lemma]{Result}
\newtheorem*{thm*}{Theorem}
\newtheorem*{prop*}{Proposition}
\newtheorem*{cor*}{Corollary}
\newtheorem*{conj*}{Conjecture}
%\numberwithin{equation}{section}
\theoremstyle{remark}
\newtheorem{rem}[lemma]{Remark}
\newtheorem*{rem*}{Remark}
\theoremstyle{definition}
\newtheorem{Def}[lemma]{Definition}
\newtheorem*{Def*}{Definition}
\def\signature#1#2{\par\noindent#1\dotfill\null\\*
{\raggedleft #2\par}}

\newcommand\eff{\mathrm{eff}}
\renewcommand{\theenumi}{\roman{enumi}}
\renewcommand{\labelenumi}{(\theenumi)}

\title{Inverse problems in $N$-body scattering}
\author[Gunther Uhlmann]{Gunther Uhlmann}
\address{Department of Mathematics, University of Washington, 
Seattle, WA}
\email{gunther@math.washington.edu}
\author[Andras Vasy]{Andr\'as Vasy}
\address{Department of Mathematics, Massachusetts Institute of Technology,
Cambridge MA 02139}
\email{andras@math.mit.edu}
\thanks{G.\ U.\ is partially supported by NSF grant \#DMS-00-70488 and
a John Simon Guggenheim fellowship.
A.\ V.\ is partially supported by NSF grant \#DMS-0201092
and a Fellowship from the Alfred P.\ Sloan Foundation.}

\maketitle

Scattering theory is the analysis of the motion of several interacting
particles.
Inverse problems in scattering theory seek the answer to the question:
can one determine the interactions between particles by a scattering
experiment? That is, if one shoots a number of particles at each other
in an accelarator and observes the outcome, can one find out how
the particles interact? In this note we attempt to explain this problem
and the existing results, and we also sketch the proof of the extension
of our recent three-body result to the many-body case:
this requires only minor modifications. We also indicate how to extend
the inverse result of \cite{Vasy:Structure} to the many-body case.

Before stating many-body results, we recall the simpler two-body setting,
where (after removal of the center of mass) one studies the Hamiltonian
$\Delta+V$
on $\Real^n$ for some $n$, where $\Delta$ is the positive Laplacian on
$\Real^n$, and $V$ is a function on $\Real^n$ that decays at infinity.
Thus, $V$
can be considered a perturbation of $\Delta$, and hence a relatively simple
analysis of scattering is possible. We recall that at high energies $V$
can be considered a perturbation in an even stronger sense. Roughly, at high
energies, $V$
is not only relatively compact, but also relatively small
(compared to $\Delta$), so one
can use the Born approximation to study inverse scattering. In particular,
it is easy to recover $V$ from the high energy ($\lambda\to+\infty$)
asymptotics of the scattering
matrix $\calS(\lambda)$ --
this is the object that captures the relationship between the
initial ($t\to-\infty$) and final ($t\to+\infty$) state of the particles
of energy $\lambda$,
and which is described below.
If one only knows the scattering matrix at a fixed energy $\lambda$,
it is still possible to recover $V$, at least
for compactly supported, or indeed exponentially decaying, potentials
\cite{Sylvester-Uhlmann:Global, Novikov:Fixed, Novikov-Khenkin:D-bar,
Uhlmann-Vasy:Fixed}. For short-range potentials that are symbols, it is
possible to recover the {\em asymptotics} of the potential at infinity
using fixed energy information only, see \cite{Joshi-SB:Asymptotics}.
For Schwartz potentials $V$ (whose asymptotics is thus trivial)
there are counterexamples
to the fixed energy problem \cite{Grinevich-Novikov:Transparent}, namely
there are `transparent' potentials,
but it would be interesting to know whether
the S-matrix given in a finite interval determines the potential, since
this would have direct applications in many-body scattering, see
Theorem~\ref{thm:3-3}.

In fact, there are many inverse problems that one can study. If there
are more than two particles, they may have a number of initial (or incoming)
and final (or outgoing) configurations.
For example, one can shoot $N$ electrons at each other, let them interact,
and observe the outcome. Since the electrons repel each other,
we will have $N$ (asymptotically) free particles as the outcome.
This is $N$-cluster to $N$-cluster scattering: there were $N$
clusters of particles (each consisting of a single electron) both
initially and finally.

In a physically much more relevant way, one may shoot
an electron at a target which is a composite particle, for example
an atom or ion. If the energy of the electron
is high enough, the ion may break up e.g.\ a helium ion may break up
into the nucleus and an electron. If the atom/ion breaks up into $k$ parts,
this is two-cluster to $k+1$-cluster scattering:
initially there were two clusters, namely the electron and the ion,
and finally there are $k+1$: the $k$ parts of the atom/ion and the initial
electron. For the helium ion, if the ion breaks up into the nucleus and an
electron, $k=2$, so this is two-cluster to three-cluster scattering.
Or, the atom/ion may not break up, though it may become excited, and
eventually one will have an atom/ion and an electron, just as initially.
This is two-cluster to two-cluster scattering. If in addition the ion
is in the same final and initial states, the problem may be instead considered
as a two-body problem, with the ion and the electron comprising the
two (now composite) particles; indeed, this is the usual way one thinks
of such an experiment physically, at least at low energies. A natural
question is then whether this is an accurate description, i.e.\ whether
the more complete three-particle analysis gives the same predictions
as the simplified two-particle approach.

Another classification of inverse problems is by the energy range
in which the scattering data are known.
High energy problems have been extensively studied by
Enss and Weder
\cite{Enss-Weder:Geometrical, Enss-Weder:Inverse},
Novikov \cite{Novikov:N-body} and Wang \cite{Wang:High, Wang:Uniqueness}.
Namely, as mentioned above, at high energies the potential
(the interaction) is small compared to the Laplacian (the kinetic energy),
in an appropriate sense. Thus, one can use the Born approximation to
study inverse problems and recover information about the interactions.
In this note we
discuss fixed energy problems, and problems where a component of the
scattering
matrix is known in a bounded interval of energies.

We first recall that the actual $N$-body Hamiltonian for $d$-dimensional
particles takes the
following form:
\begin{equation}\label{eq:mb-Ham}
H=\sum_{i=1}^N\frac{\hbar^2}{2m_i}\Delta_{q_i}+\sum_{i<j} V_{ij}(q_i-q_j),
\end{equation}
where $q_i\in\Real^d$ is the position and $m_i$ is the mass of particle $i$,
while $V_{ij}$ is the interaction between particles $i$ and $j$,
and $\hbar$ is Planck's constant.
Note that $H$ is an operator on (functions on) $\Real^{Nd}$. The quantum
$N$-body problem is then the analysis of the Schr\"odinger equation
$i\hbar \pa_t u=Hu$, i.e.\ of the propagator $e^{-iHt}$, as $t\to\pm\infty$,
or equivalently, by taking the $t$-Fourier transform, of tempered
distributional solutions of $(H-\lambda)u=0$.

It is convenient to generalize the framework, introducing some notation due to
Agmon and Sigal; see \cite{Derezinski-Gerard:Scattering} for a detailed
discussion of the setup. We work on a vector space $X_0$ endowed with a
translation invariant metric $g$; by introduction of an orthonormal basis
we may of course identify it with $\Real^n$, as we usually do below.
Let $\Delta\geq 0$ be the positive
Laplacian of $g$.
We are also given a finite collection $\calX=\{X_a:\ a\in I\}$
of linear subspaces $X_a$ of $X_0=\Rn$, called the collision planes.
It is convenient to assume that $\calX$ is closed under intersections,
and includes $X_0=\Real^n$ and $X_1=\{0\}$. We
let $X^a=X_a^\perp$ be the orthocomplement of $X_a$ in $\Rn$, so
$\Rn=X_a\oplus X^a$. We write the corresponding coordinates as $(x_a,x^a)$,
and denote the projection to $X^a$ by $\pi^a$. A many-body
Hamiltonian is an operator of the form
\begin{equation*}
H=\Delta+\sum_a (\pi^a)^*V_a,
\end{equation*}
where $V_a$ is a real valued function on $X^a$ in a certain class,
for example $V_a$ is a symbol on $X^a$ of negative order:
$V_a\in S^{-\rho}(X^a)$, $\rho>0$. We drop the pull-back notation
from now on and write $H=\Delta+\sum_a V_a$.
Since here we are interested in spectral theory, and $V_0$ is a function
on a point, hence a constant, it is convenient to assume that $V_0=0$:
otherwise it would simply shift the spectral parameter.

In the standard example, \eqref{eq:mb-Ham},
the collison planes are $X_{ij}=\{q_i=q_j\}$, as well as their intersections.
Thus,
at $X_{ij}$ particles $i$ and $j$ are at the same place, hence the name
`collision plane'.

The main feature of many-body problems is that even if $V_a$ decays at infinity
on $X^a$, it does {\em not} decay at infinity in $\Rn$ since it is
a constant along $X_a$ as well as along its translates $X_a+\{\bar x^a\}$,
$\bar x^a\in X^a$,
so it does not decay if we go to infinity, say, along $X_a$.

In the two-body problem one actually has $H=\Delta_{q_1,q_2}+V_{12}(q_1-q_2)$,
i.e.\ $V_{12}$ still does not decay at infinity, e.g.\ if one keeps
$q_1=q_2$ but lets $q_1\to\infty$. However, one can easily remove the
center of mass by performing a Fourier transform along $X_{12}$,
and then one reduces the problem to the study of a Hamiltonian on $X^{12}$
of the form
$\Delta+V_{12}$ where $V_{12}$ decays at infinity (we are working
on $X^{12}$!). This reduced Hamiltonian can be considered as
a perturbation of $\Delta$, hence its analysis is rather simple.

The center of mass can also be removed in any actual many-body problem,
but one still obtains a Hamiltonian with non-decaying potentials as before.

We also need to introduce subsystem Hamiltonians to describe bound states
such as atoms and ions. We order the clusters by inclusion
of the $X^a$, so
\begin{equation*}
X^a\subset X^b\Miff a\leq b.
\end{equation*}
For the collision planes the relation is thus reversed,
so $X_a\subset X_b$ if and only if $a\geq b$. We say that $a$ is a $k$-cluster
if the maximum number of elements in a chain with minimal element $a$ is
$k$. Thus, if $X_a$ is such that $X_b\subset X_a$ implies $b=a$ or $b=1$,
then $a$ is a 2-cluster, while $1$ is a 1-cluster.
If $0$ is an $N$-cluster, we call $H$ an $N$-body
Hamiltonian. It is easy to see that this agrees with the usual terminology
for the physical $N$-body Hamiltonian \eqref{eq:mb-Ham}, with
center of mass removed. For example, in a three-body problem one gets
that $X_a\cap X_b=\{0\}$
if $a\neq b$, $a,b\neq 0,1$.

Corresponding to each cluster $a$ we introduce the cluster
Hamiltonian $H^a$ as an operator on $L^2(X^a)$ given by
\begin{equation}
H^a=\Delta_{X^a}+\sum_{b\leq a} V_b=\Delta_{X^a}+V^a,
\end{equation}
$\Delta_{X^a}$ being the Laplacian of the induced metric on $X^a$.
Thus, if $H$ is a $N$-body Hamiltonian and $a$ is a $k$-cluster,
then $H^a$ is a $(N+1-k)$-body Hamiltonian. The $L^2$
eigenfunctions
of $H^a$ (also called bound states)
play an important role in many-body scattering.
We denote $L^2$ eigenvalues and normalized eigenfunctions of $H^a$
by $\ep_\alpha$ and $\psi_\alpha$, so $(H^a-\ep_\alpha)\psi_\alpha=0$.
The pair $(a,\psi_\alpha)$ is called a channel.
We
remark that by a result of Froese and Herbst, \cite{FroExp},
$\ep_\alpha\leq 0$ for all $\alpha$ (there are no positive
eigenvalues). Moreover, $\pspec(H^a)$ is bounded below since
$H^a$ differs from $\Delta$ by a bounded operator. Note that
$X^0=\{0\}$, $H^0=0$, so the unique eigenvalue of $H^0$ is $0$.

The eigenvalues of $H^b$ can be used to define the set of
thresholds of $H^a$. Namely, we let
\begin{equation}
\Lambda_a=\cup_{b<a}\pspec(H^b)
\end{equation}
be the set of thresholds of $H^a$, and we also let
\begin{equation}
\Lambda'_a=\Lambda_a\cup\pspec(H^a)
=\cup_{b\leq a}\pspec(H^b).
\end{equation}
Thus, $0\in\Lambda_a$ for $a\neq 0$ and $\Lambda_a\subset(-\infty,0]$.
It follows from the Mourre
theory (see e.g.\ \cite{FroMourre, Perry-Sigal-Simon:Spectral})
that $\Lambda_a$ is closed,
countable, and $\pspec(H^a)$ can only accumulate at $\Lambda
_a$.
In addition, $L^2$ eigenfunctions
of $H^a$ with eigenvalues which are not thresholds are necessarily
Schwartz functions on $X^a$ (in fact, they decay exponentially,
see \cite{FroExp}).
We write
\begin{equation*}
I_a=\sum_{b\not\leq a}V_b
\end{equation*}
for the intercluster interaction. Thus
$I_a$ is a function on $\Real^n$.

The point of intruducing $H^a$ and $H_a=\Delta_{X_a}+H^a$ is that near
\begin{equation*}
X_{a,\reg}=X_a\setminus X_{a,\sing},\ \text{where}\ X_{a,\sing}=\cup_{b>a} X_b,
\end{equation*}
$H_a$ is a good model for $H$ since $I_a=H-H_a$ decays at infinity
along $X_{a,\reg}$.
(Here we are
using that if $b\not\leq a$, then the collision plane $X_c=X_b\cap X_a$
satisfies $c>a$, i.e.\ $X_c\subsetneq X_a$.)

If we consider the two-cluster to two-cluster problem as a two-body problem by
regarding the atom/ion as a single particle, we essentially work on $X_a$
by projecting to a bound state $\psi_\alpha$ of $H^a$ and using the
effective interaction
\begin{equation*}
V_{\alpha,\eff}=\int_{X^a} I_a |\psi_\alpha|^2\,dx^a,
\end{equation*}
which is a function on $X_a$. For various reasons, including the results
in high energy settings, it is clear that we can only hope to recover
$V_{\alpha,\eff}$ (rather than $I_a$) by studying the corresponding
(i.e.\ state $\alpha$ to itself) two-cluster to two-cluster
scattering matrix.

We can now describe the S-matrices following \cite{Vasy:Scattering}.
We warn the reader that the normalization this gives is geometric
as in \cite{RBMSpec}, and
is different from the standard wave-operator normalization.
The difference is essentially given
by pull-back by an antipodal map on the sphere.
Thus, free particles propagate in
straight lines, i.e.\ they leave in the opposite direction from where they
came, so in the geometric formulation, in the absence of interactions,
the outgoing data are supported in the image of the support of the
incoming data under the antipodal map. An incoming wave
of energy $\lambda$ in the channel $\alpha$ is a function of the form
\begin{equation*}
u_{\alpha,-}=e^{-i\sqrt{\lambda-\ep_\alpha}|x_a|}|x_a|^{-\frac{\dim X_a-1}{2}}
g_{\alpha,-}(\frac{x_a}{|x_a|})\psi_\alpha(x^a)+\text{faster decaying terms},
\end{equation*}
and an outgoing wave has the form
\begin{equation*}
u_{\alpha,+}=e^{i\sqrt{\lambda-\ep_\alpha}|x_a|}|x_a|^{-\frac{\dim X_a-1}{2}}
g_{\alpha,+}(\frac{x_a}{|x_a|})\psi_\alpha(x^a)+\text{faster decaying terms},
\end{equation*}
i.e.\ the sign of the phase has changed. Here $g_{\alpha,\pm}$ may
be taken e.g.\ $L^2$ functions on $\sphere_a$, the unit sphere in $X_a$, or
ideally, at least one of them may be taken $\Cinf$.
One is then interested in generalized eigenfunctions
of $H$, i.e.\ tempered distributions $u$ on $\Rn$ that solves
$(H-\lambda)u=0$. A typical example is of the form
\begin{equation}\label{eq:gen-ef}
u=u_{\alpha,-}-(H-(\lambda+i0))^{-1}((H-\lambda)u_{\alpha,-});
\end{equation}
here the faster decaying terms may even be dropped without affecting $u$
and $g_{\alpha,-}$ can be specified to be any smooth function on $\sphere_a$.
In general, even if the incoming data are in a single channel $\alpha$,
as in \eqref{eq:gen-ef},
the corresponding generalized eigenfunction $u$ of $H$ will have outgoing waves
in all channels. The S-matrix $\calS_{\alpha\beta}(\lambda)$ picks out
the component in channel $\beta$ by projection in a certain sense,
see \cite{Vasy:Scattering}.
Thus, $\calS_{\alpha\beta}(\lambda)$ maps functions on $\sphere_a$, the unit
sphere in $X_a$, to functions on $\sphere_b$, by
\begin{equation*}
\calS_{\alpha\beta}(\lambda)g_{\alpha,-}=g_{\beta,+}
\end{equation*}
for $u$ as in \eqref{eq:gen-ef}.
For example, the free-to-free (i.e.\ $N$-cluster to $N$-cluster in
$N$-body scattering) S-matrix $\calS_{00}(\lambda)$ maps functions on $\sphere_0$,
the unit sphere in $\Rn$, to functions on $\sphere_0$, more precisely
$\calS_{00}(\lambda):L^2(\sphere_0)\to L^2(\sphere_0)$ is bounded.
Below we write $\sphere_{a,\reg}=\sphere_a\cap X_{a,\reg}$.

The first many-body result we state, slightly informally,
is an `old' direct result of the second
author \cite{Vasy:Structure} and Hassell \cite{Hassell:Plane}, which in turn
can be considered as the extension of a result of Melrose and Zworski
\cite{RBMZw} on asymptotically Euclidean spaces (i.e.\ in geometric two-body
scattering).
Its essential ingredient is the propagation of singularities, here meaning
microlocal lack of decay at infinity, much as for real principal type
PDEs in the usual microlocal setting, shown in the work of H\"ormander
\cite{Hormander:Existence, Hor}
and Kashiwara-Kawai \cite{Kashiwara-Kawai:Microhyperbolic}
for $\Cinf$ resp.\ analytic singularities
in the absence of boundaries, by
Melrose and Sj\"ostrand \cite{Melrose-Sjostrand:I, Melrose-Sjostrand:II}
and Taylor \cite{Taylor:Grazing} for smooth boundaries and $\Cinf$
singularities,
and for analytic singularities by Sj\"ostrand \cite{Sjostrand:Propagation-I},
and by Lebeau~\cite{Lebeau:Propagation} for analytic singularities at
corners. The relationship is that two-body scattering is a simpler
version of the boundariless setting \cite{RBMSpec, RBMZw},
as shown by the Fourier transform,
while 3-body scattering is analogous to the problems with smooth boundaries,
and more than three bodies correspond to the corners.

\begin{thm}\label{thm:3-3}
Suppose that $H$ is a three-body Hamiltonian and
the $V_a$ are Schwartz on $X^a$ for all $a$. Then
$\calS_{00}(\lambda)$ is a finite sum of Fourier integral operators (FIOs)
associated to the broken geodesic relation on $\sphere_0$ to distance $\pi$.
Its canonical relation corresponds to the various collision patterns.
The principal symbol of the term corresponding to a single collision at
$X_a$ is given by, and in turn determines,
the 2-body S-matrix of $H^a$ at energies $\lambda'
\in(0,\lambda)$.
\end{thm}

This result presumably extends to short range symbolic potentials, using the
same methods, though it is technically more complicated to write down the
argument in that case, and it has not been done.
However, it was shown in \cite{Vasy:Propagation-Many} and
\cite{Vasy:Bound-States} that in $N$-body scattering, $N$ arbitrary,
the wave front relation of all S-matrices
is given by the broken bicharacteristic relation. Thus, the {\em location}
of the singularities of the Schwartz kernel of the S-matrix is known,
but their {\em precise form,} e.g.\ that they are those of an FIO,
is not known except in three-body scattering. (Also, in {\em some} parts
of the wave front relation the singularities can be described precisely
in the $N$-body setting, as we do below in Theorem~\ref{thm:N-N}.)

An immediate corollary, when combined with two-body results
(e.g.\ analyticity of the S-matrix in $\lambda'$ and the Born approximation)
is the following inverse result.

\begin{cor}
In three-body scattering,
if the $V_a$ decay exponentially and $\dim X_a\geq 2$ for all $a$
then $\calS_{00}(\lambda)$ for a single value of $\lambda$ determines all
interactions.
\end{cor}

This result is analogous to the recovery of cracks in a material
by directing sound waves at it and observing the singularities
of the reflected waves, except the
last step which uses two-body results to get the potentials from the
two-body S-matrices -- and apart from the fact that the crack recovery
result does not seem to exist in the literature!
Note
that in 3-body scattering (and also in $N$-body scattering if the
set of thresholds is discrete)
there is only a finite number of possible collision patterns.
Generalized
broken bicharacteristics are curves in a compressed phase space described
in \cite{Vasy:Bound-States}, or in a more leasurely way in
\cite{Vasy:Geometry}. A subset of this family
is the set of broken bicharacteristics, which are piecewise bicharacteristics
of the metric function $g$, i.e.\ integral curves of its Hamilton vector
field in $T^*X_0=(X_0)_x\times (X_0^*)_\xi$,
that satisfy in addition that at each break point at $X_{b,\reg}\times X_0$,
the momenta $\xi$, $\xi'$ of the two segments meeting there differ by
an element of $(X^b)^*$, i.e.\ $\xi-\xi'$ vanishes on (or is conormal to)
$X_b$. This means that in a
collision at $X_{b,\reg}$,
the energy and the component of the momentum tangential to
$X_b$ are conserved, similarly to geometric optics. A convenient way of
making the broken bicharacteristics continuous is to compress the phase
space, by identifying points $(x,\xi)$ and $(x,\xi')$ if $x\in X_{b,\reg}$
and $\xi-\xi'\in (X^b)^*$. Such a description is in fact necessary
for generalized broken bicharacteristics, but the new analysis below
involves only broken bicharacteristics, so we adopt the more naive piecewise
continuous curve approach and work in $X_0\times X_0^*$.
(In 3-body scattering every generalized broken bicharacteristic is
such a broken bicharacteristic, but this is not true in general.) Note
that a bicharacteristic segment of $g$ is a curve of the form
\begin{equation}\label{eq:free-bich}
x(t)=\bar x+2t\bar\xi,\ \xi=\bar\xi,\ t\in J,\ J\ \text{an interval},
\end{equation}
i.e.\ is a straight line segment with constant momentum. In addition,
$|\xi|^2$ is the kinetic energy, so in the free region, $X_{0,\reg}$,
$|\xi|^2=\lambda$, but for particles in a bound state of energy
$\ep_\alpha$, we actually
have $|\xi|^2+\ep_\alpha=\lambda$.

\begin{figure}
\begin{center}
\mbox{\epsfig{file=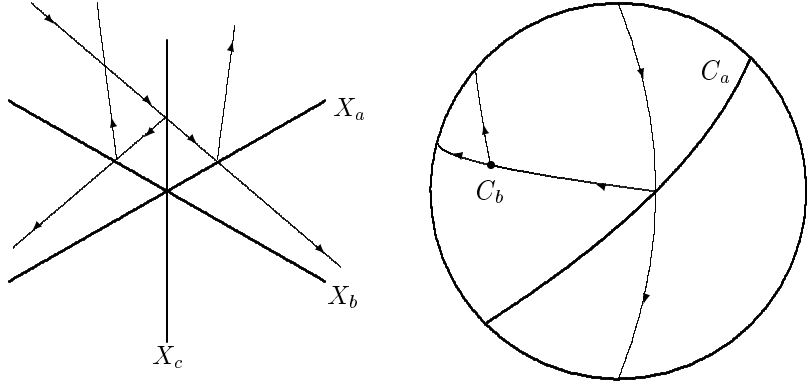}}
\end{center}
\caption{On the left,
broken geodesics in $\Real^n\setminus\{0\}$, $n=2$,
broken at the collision planes
$X_a$, $X_b$ and $X_c$. On the right, the projection of broken
geodesics in $\Real^n\setminus\{0\}$, $n=3$, emanating from the north pole,
to the unit sphere $\sphere_0$, better understood
as the sphere at infinity. The $C_a$, $C_b$ are the intersection
of the collision planes $X_a$, $X_b$ with $\sphere_0$; $\dim X_a=2$, $\dim X_b=1$.}
\label{fig:brtraj}
\end{figure}

In fact, in \cite{Vasy:Bound-States} the phase space is compactified, much
as in \cite{RBMSpec}, which here essentially means a quotient by the
$\Real^+$-action on $X_0\times X_0^*$ that is simply the dilation in the $X_0$
component. Taking $\sphere_0\times X_0^*$ as a transversal
to this action, the bicharacteristic segments become curves on $\sphere_0\times
X_0^*$. {\em Below we freely use this identification without further commenting
on it.} Then the projection of a broken bicharacteristic segment to
$\sphere_0$ is a (reparameterized) geodesic, and in
3-body scattering the total length of these segments for a maximally extended
broken bicharacteristic is $\pi$, see Figure~\ref{fig:brtraj}.
(One example, corresponding to an unbroken
bicharacteristic, i.e.\ a straight line, is half of a great circle on
$\sphere_0$, connecting antipodal points.)

Unfortunately, this was necessarily a very brief and incomplete
discussion, and we refer
the reader to \cite{Vasy:Bound-States, Vasy:Geometry} for a much more detailed
description of generalized broken bicharacteristics.

While the $N$-cluster to $N$-cluster S-matrix is more complicated
in general than for $N=3$, and in particular the FIO result is not
known, it was proved in \cite{Vasy:Bound-States} that the wave
front relation of the S-matrices is given by the generalized broken
bicharacteristic relation. These generalized broken bicharacteristics
can be quite complicated in the presence of bound states; not all of
them are necessarily broken bicharacteristics. However, it turns out
that for an open dense subset of
the relation corresponding to a single collision
in an $(N-1)$-cluster $b$
there are {\em no} generalized broken bicharacteristics relating these points
except those with a
single break at $b$.

Indeed, it is easy to write
down this relation, $\Lambda_b$; this was done, for example, in
\cite[Section~6]{Vasy:Structure}. Recall that when discussing canonical
relations, such as $\Lambda_b\subset T^*\sphere_0\times T^*\sphere_0$,
it is convenient to use the twisted symplectic form $\omega-\omega'$
where $\omega$, resp.\ $\omega'$, is the standard symplectic form
on the first, resp.\ second, copy of $T^*\sphere_0$.
We identify points in $(y,\mu)\in T^*\sphere_0$
with pairs $(y,\mu)\in\sphere_0\times X_0=\sphere_0\times\Rn$ such that
$y\cdot\mu=0$, and write
$P^\perp_y$ for the orthogonal projection to the orthocomplement
to $\Span\{y\}$. Then the relation
is the twisted conormal bundle of
\begin{equation}\label{eq:refl-pts}
\{(y,y'):\ y'_b=-y_b,\ y'\in\sphere_0
\setminus\sphere_b\}
\end{equation}
with the zero section removed
(twisting means
that the sign of the dual variable of $y'$ is changed); the
restriction $y'\in\sphere_0
\setminus\sphere_b$ ensures that \eqref{eq:refl-pts} is a smooth
manifold. Its twisted conormal bundle is
\begin{equation}\begin{split}\label{eq:Lambda_b-0-def}
\Lambda_b=\{(y,\mu,y',\mu'):\ &y'\in\sphere_0\setminus\sphere_b,
\ y_b=-y'_b,\ \mu_b=-\mu'_b,\\
&\mu^b=Cy^b,\ (\mu')^b=-C(y')^b,\ C=-\frac{\mu_b\cdot y_b}{|y_b|^2}\},
\end{split}\end{equation}
i.e.\ is parameterized by $(y', y^b, \mu_b)$ with
$y'\in\sphere_0\setminus\sphere_b$, $y^b\in X^b$ satisfying
$|y^b|^2=1-|y'_b|^2$, and $\mu'_b\in X_b\setminus\{0\}$.
Note that $C$ is determined
by $y\cdot\mu=0=y'\cdot\mu'$.

The (unbroken) geodesic of length $\pi$
on $\sphere_0$ associated to $(y',\mu')$ is the arc given by $(\Span\{y'\}+
\{c\mu':\ c>0\})\cap\sphere_0$; it is thus the great circle starting at $y'$
in the direction $\mu'$. Rather than using $\mu'_b$ as a parameter in
\eqref{eq:Lambda_b-0-def}, it is natural to use $z\in\sphere_b$ and
$|\mu'|=|\mu|$, where $z$ is the point at which the broken geodesic hits
$\sphere_b$. Thus, $\mu'=cP^\perp_{y'}z$, $c>0$, and an
alternative description of $\Lambda_b$ is
\begin{equation}\begin{split}\label{eq:Lambda_b-def}
\Lambda_b=\{(y,\mu,y',\mu'):\ &\exists z\in \sphere_b,\ c>0,
\ y'\in\sphere_0\setminus\sphere_b,\ y_b=-y'_b,\\
&\mu=-cP^\perp_y z,\ \mu'=cP^\perp_{y'}z\}.
\end{split}\end{equation}
Thus, $y\in\sphere_{0,\reg}$, $z\in\sphere_b$, $y^b\in X^b$ with
$|y^b|^2=1-|y'_b|^2$, and $|\mu|$ parameterize $\Lambda_b$.
(It is straightforward to express $z$ in terms of $(y',\mu')$ as well,
and show that the change of variables $(y',\mu')\to(y',z)$ is smooth.)
The reader may find it helpful to study the parameterizations
\eqref{eq:bich-out-param}-\eqref{eq:bich-in-param} below.

We also $\Lambda_0$ be the free bicharacteristic
relation, which is thus the twisted conormal bundle of $y=-y'$.

Let
\begin{equation}\begin{split}\label{eq:Lambda_b-p-def}
\Lambda'_b=
\{(y,\mu,y',\mu'):\ \exists z\in\sphere_{b,\reg},&\ c>0,\ y_b=-y'_b,
\ \mu=-cP^\perp_y z,\ \mu'=cP^\perp_{y'}z,\\
&y\neq -y',
\ y,y'\nin\Span(z)\oplus
X_a\ \text{for any}\ a\neq 0,b\}.
\end{split}\end{equation}
Then we have the following
new result, generalizing Theorem~\ref{thm:3-3}.

\begin{thm}\label{thm:N-N}
Suppose that $H$ is an $N$-body Hamiltonian and
the $V_a$ are Schwartz on $X^a$ for all $a$, $\dim X^a\geq 2$ for all $a$.
With $b$ an $(N-1)$-cluster, $\Lambda_b$, $\Lambda'_b$ as above,
$\Lambda'_b$ is an open dense subset of $\Lambda_b$. Moreover,
$\calS_{00}(\lambda)$
is a Fourier integral operator microlocally at $\Lambda'_b$, and its
principal symbol at $(y,\mu,y',\mu')\in\Lambda_b$ is determined by, and in turn
determines, the Schwartz kernel of the
S-matrix $\calS^b(\lambda|y^b|)$
of $H^b$ (which is a 2-body Hamiltonian)
at $(\frac{y^b}{|y^b|},\frac{(y')^b}{|y^b|})$.
In particular, $\calS_{00}(\lambda)$ determines
$\calS^b(\sigma)$ for $\sigma\in(0,\lambda]$, hence also $V_b$, provided that
$V_b$ decays exponentially.
\end{thm}

\begin{proof}
First, openness is clear since each additional condition in the definition
\eqref{eq:Lambda_b-p-def} of
$\Lambda'_b$ (as compared to \eqref{eq:Lambda_b-def})
gives an open subset of $\Lambda_b$, and the intersection of these is open.

To see density,
for $z\in \sphere_{b,\reg}$, note that $X_a\oplus\Span(z)$ has codimension
$\dim X^a-1\geq 1$ in $X_0$, so imposing the additional conditions
$z\in\sphere_{b,\reg}$ and
$y'\nin\Span(z)\oplus
X_a$ for any $a\neq 0,b$
certainly gives a dense subset of $\Lambda_b$.

Consider $y'\in\sphere_{0,\reg}
\setminus \cup_{a\neq b,0}(X_a\oplus \Span(z))$. Then a point in
$\Lambda_b$ is specified by giving $y\in\sphere_0$ satisfying
$y_b=-y'_b$, and $|\mu|$. The projection of the
corresponding bicharacteristic segments to the base space, $\sphere_0$,
i.e.\ the corresponding geodesics connecting $y'$ and $z$, resp.\ $y$ and $z$,
lie in $\Span\{y',z\}\cap\sphere_0$ and $\Span\{y,z\}\cap\sphere_0$.
If these intersect
$X_a$ then $sy'+tz\in X_a$, resp.\ $sy+tz\in X_a$, for some $(s,t)\in\Real^2$
so
$y'\in X_a\oplus\Span\{z\}$, resp.\ $y\in X_a\oplus\Span\{z\}$,
as $s\neq 0$. (For $s=0$ gives $z\in X_a$, contradicting
$z\in \sphere_{b,\reg}$.) Note that $y'\nin X_a\oplus\Span\{z\}$ by assumption.
It remains to show that the set of $y\in\sphere_{0,\reg}$ with $y_b=-y'_b$ and
$y\nin X_a\oplus\Span\{z\}$ is also dense inside
\begin{equation*}
S=\{y\in\sphere_0:\ y_b=-y'_b\}.
\end{equation*}
Thus, $S$ is a sphere inside the affine space
$\{-y'_b\}+X^b$, centered
at $-y'_b$, of radius $\sqrt{1-|y'_b|^2}$. The other conditions on $y$
state that it does not lie on various affine subspaces $W_a$, namely
$X_a\oplus\Span\{z\}$, $a\neq 0,b$, and $W_b=X_b$,
and in addition $y\neq -y'$.
But the intersection
of a sphere with an affine space $W$ is a lower dimensional sphere, unless
the whole sphere lies in $W$. Since $y=-y'$
does not lie on any of these affine subspaces $W_a$, $a\neq 0$,
all of these intersections
are lower dimensional spheres, so in fact the set of $y$
satisfying $y\nin X_a\oplus\Span\{z\}$ and $y\in\sphere_{0,\reg}$ is
the complement of the union of compact codimension $\geq 1$ submanifolds
of $S$, and hence is dense in $S$. (A slightly different way of arguing
is to say that if a point $\bar y$ lies in one of these affine spaces, then
the arc between $\bar y$ and $-y'$ only meets this space inside
$\Span\{\bar y\}$, for otherwise $-y'$ would also lie in this affine space,
giving the density rather explicitly.) Since $y=-y'$ is a single point in $S$,
and $S$ has dimension $\dim X^b-1\geq 1$,
the density follows even if we impose $y\neq -y'$.

The bicharacteristics \eqref{eq:free-bich}
of the free Hamiltonian $\Delta$, at energy $\lambda$,
can be parameterized by $p=(y,\mu)\in T^*\sphere_0$
by
\begin{equation}\label{eq:bich-out-param}
\gamma_p(t)=(x(t),\xi(t)),\ x(t)=-\frac{\mu}{\sqrt{\lambda}}+2t\xi,
\ \xi(t)=\sqrt{\lambda}\,y.
\end{equation}
The analogous `incoming' parameterization, $\gamma'_{p'}(t)=(x'(t),\xi'(t))$,
$p'=(y',\mu')$,
simply switches the sign of $y$
and $\mu$, i.e.\ it is
\begin{equation}\label{eq:bich-in-param}
x(t)=\frac{\mu'}{\sqrt{\lambda}}+2t\xi,\ \xi(t)=-\sqrt{\lambda}\,y',
\end{equation}
so as $t\to+\infty$, $\frac{x(t)}{|x(t)|}\to y$ along $\gamma_p$,
while as $t\to-\infty$, $\frac{x(t)}{|x(t)|}\to y'$ along $\gamma'_{p'}$.
The normalization of these
parameterizations may seem strange; they arise by considering the
graph of $d(\sqrt{\lambda}x\cdot y)$, resp.\ $d(-\sqrt{\lambda}x\cdot y')$,
(recall that $e^{\pm i\sqrt{\lambda}x\cdot y}$ is, up to constants, the
Schwartz kernel of the free Poisson operator!), expressing $x$ and $\xi$
in terms of $(y,\mu)$ and an additional parameter $t$, and then changing
the sign of the $\mu$ component, corresponding to the twisting.

Now suppose that $(p,p')\in\Lambda'_b$, $p=(y,\mu)$, $p'=(y',\mu')$.
Then there is a unique
generalized broken bicharacteristic `connecting' $p$ and $p'$, i.e.\ with
$\gamma$ such that for large $|t|$, up to reparameterization, $\gamma$
is given by $\gamma_p$ and $\gamma'_{p'}$; say
$\gamma|_{[T,+\infty)}$ coincides with $\gamma_p$
and $\gamma|_{(-\infty,T']}$ coincides with $\gamma'_{p'}$, up to
reparameterization.
Moreover, a generalized broken bicharacteristic
(being a bicharacteristic in $T^*X_{0,\reg}$) has a unique continuation
until it hits $\sphere_c$ for some $c\neq 0$, so we may assume
that $\gamma(T),\gamma(T')\in\sphere_{0,\sing}$.
Note that $\gamma_0$, the broken bicharacteristic
with a single break at $\sphere_b$ (with some parameterization),
is such a bicharacteristic, so $\gamma(t)=\gamma_0(t-t_+)$, $t\geq T$,
and $\gamma(t)=\gamma_0(t-t_-)$, $t\leq T'$. In particular,
$\gamma_0(T'-t_-),\gamma_0(T-t_+)\in\sphere_{0,\sing}$. Since $\gamma_0$ has
a unique such point, they are equal, so $T'-t_-=T-t_+$, and
$q=\gamma(T')=\gamma(T)$. Since the radial momentum, $\frac{x\cdot\xi}{|x|}$,
is strictly
increasing along $\gamma$ away from the radial points, in particular
at $q$, we deduce that $T'=T$, hence $\gamma=\gamma_0$ as claimed. (Recall
that radial points are those where the Hamilton vector field is radial,
which in this case means that for $x\in X_{b,\reg}$, $\xi_b=c x$ for some
$c\in\Real$.)

To construct $\calS_{00}(\lambda)$ microlocally, we construct the Poisson
operators first. Microlocally a good approximation
for these are the Poisson operators $\calP^\sharp_{0,+}(\lambda)$ of $H_b$,
which are easy to write down explicitly, see \cite{Vasy:Structure}, since
$H_b$ is `product-type': a 2-body Hamiltonian $H^b$ with the `center
of mass' kinetic energy $\Delta_{X_b}$ added back.
For suitable microlocal cutoffs $T_+$ (which are thus ps.d.o.'s)
\begin{equation*}\begin{split}
&\calP_{0,+}(\lambda)g=\calP^\sharp_{0,+}(\lambda)g
-R(\lambda+i0)I_b\calP^\sharp_{0,+}(\lambda)g,\\
&\calS_{00}(\lambda)g=((H-\lambda)T_+\calP^\sharp_{0,-}(\lambda))^*
\calP_{0,+}(\lambda);
\end{split}\end{equation*}
see \cite{Vasy:Scattering} or \cite[Equation~(4.2)]{Vasy:Geometry}.
Now, the wave front relation of
\begin{equation*}
((H-\lambda)T_+\calP_{0,-}(\lambda))^*
R(\lambda+i0)I_b\calP^\sharp_{0,+}(\lambda)
\end{equation*}
is a subset of the
generalized broken bicharacteristic relation, but it includes only those
generalized broken bicharacteristics on which $I_b$ is not Schwartz,
i.e.\ which go through points over $\sphere_a$ with $a\neq 0,b$. (See the
proof of Theorem~2.8 in \cite[Section~10]{Vasy:Bound-States} for a similar
argument.) In particular,
it does not contain any points in $\Lambda'_b$, so microlocally at
$\Lambda'_b$, $\calS_{00}(\lambda)$ is given by
\begin{equation}\label{eq:N-N-8}
((H-\lambda)T_+\calP^\sharp_{0,-}(\lambda))^*\calP^\sharp_{0,+}(\lambda)
\end{equation}
But this is an FIO by exactly the same proof as in the three-body
setting \cite{Vasy:Structure}, proving the FIO statement.

The rest of the proof is exactly as in \cite[Proposition~6.4]{Vasy:Structure},
using the explicit form of \eqref{eq:N-N-8}.
\end{proof}

The other
new result is in 2-cluster to 2-cluster scattering. To set it
up, we first state a theorem of Skibsted \cite{Skibsted:Smoothness}.

\begin{thm*} \cite{Skibsted:Smoothness}
The two-cluster to two-cluster S-matrices have $\Cinf$ Schwartz kernels,
except the conormal singularity of $\calS_{\alpha\alpha}(\lambda)$
corresponding to free motion, i.e.\ at the graph of the antipodal map
on $\sphere_a$.
\end{thm*}

Thus, principal symbol calculations do not help in this inverse problem.
The new result is the following extension of \cite{Uhlmann-Vasy:Low} to the
$N$-body setting.

\begin{thm}\label{thm:main}
Suppose that $a$ is a 2-cluster, $\dim X_a\geq 2$,
and $V_b$ for $b\leq a$ is a
symbol of negative order (i.e.\ may be long range).
Suppose also that
\begin{equation*}
\inf\spec H^a<\inf\Lambda_a
\end{equation*}
(and hence $\inf\spec H^a$ is an $L^2$-eigenvalue of $H^a$),
$I\subset(\ep_\alpha,\inf\Lambda_a)$
is a non-empty open set with $\sup I<\inf\Lambda_a$, and let
\begin{equation*}
R=2\sqrt{\sup I-\ep_\alpha}.
\end{equation*}
For any
$\mu>\dim X_a$ there
exists $\delta>0$ such that the following holds.

Suppose that $\sup|(1+|x^b|)^\mu V_b(x^b)|<\delta$ for all $b$ with
$b\not\leq a$.
Then $\calS_{\alpha'\alpha''}(\lambda)$ given for all $\lambda\in I$ and
for all bound states $\alpha',\alpha''$ of $H^a$ with $\ep_{\alpha'},
\ep_{\alpha''}<\sup I$, determines the
Fourier transform of the effective interaction $V_{\alpha,\eff}$
in the ball of radius $R$ centered at $0$.
\end{thm}

\begin{rem*}
Here $\delta$ depends on $I$ only through $\sup I$; as long as
$\sup I$ is bounded away from $\inf\Lambda_a$, $\delta$ may be chosen
independently of $I$. This dependence on $\sup I$ was improperly
omitted from the statement of Theorem~1.1 in \cite{Uhlmann-Vasy:Low}.

The smallness assumption on the $V_b$, $b\not\leq a$ implies that
$\inf\spec H^b$ is not much lower than $\inf\Lambda_a$, in particular
it may be assumed to lie above $\sup I$.
\end{rem*}

This theorem says that if the unknown interactions are small then the
effective interaction can be partly determined from the knowledge of
all S-matrices with incoming and outgoing data in the cluster $a$ in
the relevant energy range.
In fact, near-forward information suffices as in two-body scattering,
where this was observed recently by Novikov \cite{Novikov:Determination},
see \cite{Uhlmann-Vasy:Low}.
Also, if one is willing to take small $R$ and
$\alpha$ is the ground state of $H^a$,
it suffices to know $\calS_{\alpha\alpha}(\lambda)$
to recover $\hat V_{\alpha,\eff}$. That is, one has

\begin{thm}\label{thm:main2}
Suppose that $a$ is a 2-cluster, $\dim X_a\geq 2$, and $V_b$ for $b\leq a$ is a
symbol of negative order (i.e.\ may be long range).
Suppose also that
\begin{equation*}
\ep_\alpha=\inf\spec H^a<\inf\Lambda_a.
\end{equation*}
Let $\ep_{\alpha'}$ be the second eigenvalue of $H^a$, or the bottom
of its essential spectrum if $\ep_{\alpha}$ is the only eigenvalue.
Suppose also that $I\subset(\ep_\alpha,\ep_{\alpha'})$
is a non-empty open set with $\sup I<\inf\Lambda_a$, and let
\begin{equation*}
R=2\sqrt{\sup I-\ep_\alpha}.
\end{equation*}
For any
$\mu>\dim X_a$ there
exists $\delta>0$ such that the following holds.

Suppose that $\sup|(1+|x^b|)^\mu V_b(x^b)|<\delta$ for all $b$ with
$b\not\leq a$.
Then $\calS_{\alpha\alpha}(\lambda)$ given for all $\lambda\in I$
determines the
Fourier transform of the effective interaction $V_{\alpha,\eff}$
in the ball of radius $R$ centered at $0$.
\end{thm}

In case $V_b$ decay exponentially on $X^b$ for all $b\not\leq a$, then
$V_{\alpha,\eff}$ decays exponentially on $X_a$, hence its Fourier
transform is analytic, so $V_{\alpha,\eff}$ itself can be recovered
from these S-matrices.

This result should extend to higher energies, i.e.\ $\sup I< \inf\Lambda_a$
is {\em not} expected to be essential. But it is hard to make
$R$ greater than $2\sqrt{\inf\Lambda_a-\ep_\alpha}$ even then.
The reason is that our method relies on the construction of exponential
solutions following Faddeev \cite{Faddeev:Inverse},
Calder\'on \cite{Calderon:Inverse}, Sylvester and the first author
\cite{Sylvester-Uhlmann:Global} and Novikov and
Khenkin~\cite{Novikov-Khenkin:D-bar}, but in the many-body setting.
One thus allows complex momenta $\rho\in\Cx(X_a)$, the complexification
of $X_a$. We then seek solutions $u$ of $(H-\lambda)u=0$ of the form
\begin{equation}\label{eq:u_rho-def}
u=u_\rho
=e^{i\rho\cdot x_a}(\psi_\alpha(x^a)+v),\ \rho\cdot\rho=\lambda-\ep_\alpha,
\ \rho\in \Cx(X_a),
\end{equation}
where $v$ is supposed to be `small'. Note that
$u_\rho^0=e^{i\rho\cdot x_a}\psi_\alpha(x^a)$ satisfies
$(\Delta+V^a-\lambda)u_\rho^0=0$ for $\rho$ as in \eqref{eq:u_rho-def},
so $(H-\lambda)u_\rho^0=I_a\psi_\alpha$ decays at infinity since $\psi_\alpha$
decays in cones disjoint from $X_a$ while $I_a$ decays on $X_a$.
(We are using here that $a$ is a 2-cluster for otherwise $I_a$ would only
decay on $X_{a,\reg}$!) Writing
\begin{equation*}
\rho=z\nu+\rho_\perp\Mwith |\nu|=1,
\rho_\perp\cdot\nu=0,\ \rho,\nu\ \text{real},\ z\in\Cx,
\end{equation*}
the atom/ion may
break up for $|\rho_\perp|>\sqrt{\inf\Lambda_a-\ep_\alpha}$
even if $\lambda<\inf\Lambda_a$,
i.e.\ where the conservation of energy does {\em not} allow this
for {\em real} frequencies. On the other hand, one needs such large
$\rho_\perp$ to recover $V_{\alpha,\eff}$ on larger balls.

We briefly indicate how the break-up can happen. One considers the
conjugated Hamiltonian
\begin{equation*}
P(\rho)=e^{-i\rho\cdot x_a}(H-\lambda)e^{i\rho\cdot x_a}
=\Delta+2\rho\cdot D_{x_a}+\rho\cdot\rho+V-\lambda
\end{equation*}
with $\rho\in\Cx(X_a)$ the complex frequency. The total energy of the
system in state $\alpha$ is $\rho\cdot\rho+\ep_\alpha$ (kinetic+potential
energy), and we are assuming that the total energy is $\lambda$, so
\begin{equation}\label{eq:energy}
\rho\cdot\rho+\ep_\alpha=\lambda,\ \text{see}\ \eqref{eq:u_rho-def},
\end{equation}
so
\begin{equation*}
P(\rho)=\Delta+2\rho\cdot D_{x_a}+V-\ep_\alpha.
\end{equation*}
A good model for $P(\rho)$, under the smallness assumption on $I_a$, is
\begin{equation}\label{eq:P_a-def}
P_a(\rho)=\Delta+2\rho\cdot D_{x_a}+V^a-\ep_\alpha.
\end{equation}
Taking the Fourier transform in $X_a$ makes the invertibility of
$P_a(\rho)$ into a question on the behavior of the resolvent
of $H^a=\Delta_{x^a}+V^a(x^a)$, uniformly across the spectrum.
That is,
\begin{equation*}
\Fr_{X_a} P_a(\rho)\Fr_{X_a}^{-1}=
\Delta_{X^a}+V^a-(\ep_\alpha-|\xi_a|^2-2\rho\cdot\xi_a),
\end{equation*}
acting pointwise in $\xi_a$. Its right inverse is given by
\begin{equation}\label{eq:G_a-def}
\Fr_{X_a} G_a(\rho)\Fr_{X_a}^{-1}=R^a(\ep_\alpha-|\xi_a|^2-2\rho\cdot\xi_a),
\end{equation}
with $R^a(\sigma)=(H^a-\sigma)^{-1}$, provided we show that this makes sense.
Now, we need to keep
\begin{equation}\label{eq:F_rho-def}
F_\rho(\xi_a)=
\ep_\alpha-|\xi_a|^2-2\rho\cdot\xi_a\in\Cx\setminus[\inf\Lambda_a,+\infty),
\end{equation}
for otherwise other channels also become open, i.e.\ the cluster may
break up. In other words, on the range of $F_\rho$, we want to keep
$R^a$ invertible on the range of $\Id-e_a$, which would reflect
that only the $a$-bound states form open channels.
We thus need that either $\im F_\rho(\xi_a)\neq 0$ or
$\im F_\rho(\xi_a)= 0$ but $\re F_\rho(\xi_a)<\inf\Lambda_a$.
Writing $\rho=z\nu+\rho_\perp$ as above,
\begin{equation*}
\im F_\rho(\xi_a)=-2(\im z)\nu\cdot\xi_a,
\end{equation*}
so for $\im z\neq 0$, i.e.\ for $\rho$ non-real, $\im F_\rho$ vanishes
exactly on the hyperplane $\nu\cdot\xi_a=0$. Moreover,
\begin{equation}\label{eq:re-F_rho}
\re F_\rho(\xi_a)=\ep_\alpha+\rho_\perp^2
-((\xi_a)_\perp+\rho_\perp)^2
-(\nu\cdot\xi_a)^2-2(\re z)(\nu\cdot\xi_a),
\end{equation}
so if $\im F_\rho(\xi_a)=0$, the condition
$\re F_\rho(\xi_a)<\inf\Lambda_a$ amounts to
\begin{equation*}
\ep_\alpha+\rho_\perp^2
-((\xi_a)_\perp+\rho_\perp)^2<\inf\Lambda_a.
\end{equation*}
This will {\em not} be satisfied for every $\xi_a$ unless
\begin{equation}\label{eq:rho_perp-bd}
|\rho_\perp|<\sqrt{\inf\Lambda_a-\ep_\alpha}.
\end{equation}

The most unfortunate restriction of the theorem is the smallness assumption
on the unknown interactions. While \eqref{eq:rho_perp-bd}
ensures that the 2-cluster $a$ may not break up into its subsystems,
it does not rule out the existence of different channels, e.g.\ associated
to other 2-clusters, at complex frequencies. The presence of such channels
would significantly complicate the analysis. In particular, the crucial
analyticity in $z$ would not be clear at all.

We proceed to sketch the proof of Theorem~\ref{thm:main}, which is
completely analogous to the three-body result shown in \cite{Uhlmann-Vasy:Low}.
For
\begin{equation}\label{eq:ep_1-def}
\evpr\in [|\rho_\perp|^2+\ep_\alpha,\inf\Lambda_a)\setminus
\Lambda'_a,
\end{equation}
let $e_a(\evpr)$ be the orthogonal projection
to the $L^2$ eigenfunctions of $H^a$
with eigenvalue $\leq\evpr$,
and let $E_a$ be its extension to $X_0$ via
tensoring by $\Id_{X_a}$, so
\begin{equation}
E_a(\evpr)=\sum_{\alpha'}
\Id_{X_a}\otimes(\psi_{\alpha'}\otimes\overline{\psi_{\alpha'}}).
\end{equation}
Since eigenvalues of $H^a$ can only accumulate at $\Lambda_a$,
$e_a$ is finite rank. We also let
\begin{equation}\label{eq:ep_0-def}
\evth=\inf\left(\Lambda'_a\cap(\evpr,+\infty)\right)>\evpr.
\end{equation}
{\em The particular choice of $\evpr$, provided that it is sufficiently
close to $\inf\Lambda_a$,
does not play a major role in our arguments, so we usually
simply write $e_a$ for $e_a(\evpr)$, etc.}

We restrict the region \eqref{eq:rho_perp-bd}
slightly further and work in the region
\begin{equation}
\Cx(X_a)_\alpha^\circ
=\{(z,\nu,\rho_\perp):\ \im z\neq 0,
\ |\rho_\perp|^2+\ep_\alpha\in(\ep_\alpha,\evth)\setminus\Lambda'_a\},
\end{equation}
i.e.\ we also assume that $|\rho_\perp|^2+\ep_\alpha$ is not an eigenvalue
of $H^a$. Again, we do not indicate $\evth$ explicitly in the notation.

Since the ranges of $E_a$ and $\Id-E_a$ play a rather different role
below, for $p\in\Real$ we introduce weighted spaces that reflect this:
\begin{equation}\begin{split}\label{eq:calH_p-def}
\calH_p=&(L^2_p(X_a)\otimes \Range e_a)\oplus (L^2(X_a)\otimes
\Range (\Id-e_a))\\
&\qquad\subset (L^2_p(X_a)\otimes \Range e_a)\oplus L^2(X_0),
\end{split}\end{equation}
with $e_a$ considered as a bounded operator on $L^2(X^a)$.
Thus, we allow weights on the range of $E_a$, but not on its orthocomplement.
Since $\Range e_a$ is finite dimensional,
$\calH_p$ is a Hilbert
space with norms induced on the summands by the $L^2_p(X_a)$ and
$L^2(X_0)$ norms respectively.

Returning to $P_a(\rho)$, we can define a right inverse $G_a(\rho)$
by \eqref{eq:G_a-def} if we show that it makes sense, the
only issues being the behavior of $R^a(\sigma)$ for real $\sigma$ and bounds as
$|\xi_a|\to\infty$. It is straightforward to analyze these,
as in \cite{Uhlmann-Vasy:Low}, and prove the following proposition.

\begin{prop}
Suppose that $(z,\nu,\rho_\perp)\in\Cx(X_a)_\alpha^\circ$.
The operator
\begin{equation*}
G_a(\rho)=\Frinv_{X_a}R^a(\ep_\alpha-|.|^2-2\rho\cdot .)\Fr_{X_a}
\end{equation*}
is a bounded operator
$\calH_p\to \calH_r$ for
$p>0$, $r<0$, $r<p-1$. It satisfies
\begin{equation}\begin{split}\label{eq:P_aG_a-Id}
&P_a(\rho)G_a(\rho)=\Id:\calH_p\to \calH_p,\\
&G_a(\rho)P_a(\rho)=\Id:\calH_p\to \calH_p.
\end{split}\end{equation}
It is continuous in $\rho\in\Cx(X_a)_\alpha^\circ$ and
analytic in $z\in\Cx\setminus\Real$.
Moreover, for $p>0$, $r<0$, $r<p-1$, $\rho_\perp$ fixed,
for any $C>0$, $G_a(\rho)$ is uniformly bounded in $\bop(\calH_p,\calH_r)$
in $|\im z|\geq C|\re z|$,
and $\slim_{|z|\to\infty} G_a(\rho)=0$
as an operator in $\bop(\calH_p,\calH_r)$,
provided that $|z|\to\infty$ in the region
$|\im z|\geq C|\re z|$.
\end{prop}

\begin{proof}[Sketch of proof.] [See \cite{Uhlmann-Vasy:Low} for the completely
analogous proof in the three-body setting.]
Recall that
We need to analyze $F_\rho$ and $R^a(\sigma)$. First, by the spectral
mapping theorem,
for $\re\sigma<\inf\spec H^a$,
\begin{equation*}
\|R^a(\sigma)\|_{\bop(L^2,L^2)}\leq (\inf\spec H^a-\re\sigma)^{-1}.
\end{equation*}
Now, due to the $-|\xi_a|^2$ term in\eqref{eq:F_rho-def},
for any fixed $\rho$ there exists $C>0$ such that
\begin{equation}\label{eq:re F<0}
\re F_\rho(\xi_a)<C-|\xi_a|^2/2.
\end{equation}
Thus, for $\phi\in\Cinf_c(X_a^*)$ identically $1$
on a large enough ball, $G_a\Frinv_{X_a}(1-\phi)\Fr_{X_a}$ is bounded
on $L^2$. Hence, we only need to be concerned about what happens in
a compact set in $X_a^*$.

For this we need two other bounds that hold
by the selfadjointness of $H^a$ and its spectral properties, namely
\begin{equation}\begin{split}\label{eq:spec-bounds}
&\|(\Id-e_a)R^a(\sigma)\|_{\bop(L^2,L^2)}\leq |\re\sigma-\evth|^{-1},
\ \re\sigma<\evth,\\
&\|R^a(\sigma)\|_{\bop(L^2,L^2)}\leq |\im\sigma|^{-1},\ \im\sigma\neq 0,
\end{split}\end{equation}
$\evth$ as in \eqref{eq:ep_0-def}.
It is easy to check using \eqref{eq:re-F_rho}
that
\begin{equation}\label{eq:sgn-im-F-p}
\re F_\rho(\xi_a)>\ep_\alpha+\rho_\perp^2\Rightarrow
|\im F_\rho(\xi_a)|>(\re F_\rho(\xi_a)-(\ep_\alpha+\rho_\perp^2))
\frac{|\im z|}{|\re z|}
\end{equation}
In particular, $(z,\nu,\rho_\perp)\in\Cx(X_a)_\alpha^\circ$ implies
that $\ep_\alpha+\rho_\perp^2<\inf\Lambda_a$, hence
\begin{equation*}
\Frinv_{X_a}R^a(F_\rho(.))(\Id-E_a)\Fr_{X_a}
\end{equation*}
is well defined. Indeed, either $\re F_\rho(\xi_a)>\ep_\alpha+\rho_\perp^2$,
and then $\im F_\rho(\xi_a)\neq 0$, so $R^a(F_\rho(\xi_a))$ is a bounded
operator on $L^2$, or $\re F_\rho(\xi_a)\leq\ep_\alpha+\rho_\perp^2<\inf
\Lambda_a$, so $(\Id-e_a)R^a(F_\rho(\xi_a))$ is bounded on $L^2$.

On the range of $e_{a,\ep}$, the projection
to the eigenspace with eigenvalue $\ep$, $R^a(\sigma)$ is multiplication
by $(\ep-\sigma)^{-1}$. This is a locally integrable function of $\sigma$
near $\ep$ (in $\Cx$!), so the application of $R^a(F_\rho(.))$ to
$\hat u=\Fr_{X_a}u$ is well defined, provided that at every point $\xi_a$
with $F_\rho(\xi_a)=\ep$, the differential of $F_\rho:X_a^*\to\Cx$
is surjective (with $\Cx$ considered as a 2-dimensional real manifold),
i.e.\ $F_\rho(\xi_a)=\ep$ implies that $d\re F_\rho(\xi_a)$
and $d\im F_\rho(\xi_a)$ are linearly independent.
But $d\im F_\rho$ is nonzero, and is a multiple of
$d(\nu\cdot\xi_a)$, so in view of \eqref{eq:re-F_rho},
$d\im F_\rho$ and $d\re F_\rho$ are linearly independent if and only if
$d(((\xi_a)_\perp+\rho_\perp)^2)\neq 0$, i.e.\ if and only if
$(\xi_a)_\perp\neq-\rho_\perp$.
But if $(\xi_a)_\perp=-\rho_\perp$ and $F_\rho(\xi_a)$ is real,
then $\nu\cdot\xi_a=0$, hence $F_\rho(\xi_a)=\ep_\alpha+\rho_\perp^2\nin
\Lambda_a'$ due to $(z,\nu,\rho_\perp)\in\Cx(X_a)_\alpha^\circ$.
Now let
\begin{equation*}\begin{split}
G_a(\rho)&=\Frinv_{X_a}R^a(\ep_\alpha-|.|^2-2\rho\cdot .)\Fr_{X_a}\\
&=\Frinv_{X_a}R^a(\ep_\alpha-|.|^2-2\rho\cdot .)(\Id-E_a)\Fr_{X_a}
+\Frinv_{X_a}R^a(\ep_\alpha-|.|^2-2\rho\cdot .)E_a\Fr_{X_a};
\end{split}\end{equation*}
both terms are well defined by the preceeding considerations when applied
to functions in $\Sch(X_a;L^2(X^a))$. Indeed, application of $G_a(\rho)$
to the range of $E_{a,\ep}$ is the only issue, and there,
with $u=v\otimes \psi_{\alpha'}$,
$v\in L^2_p(X_a)$, $\hat v=\Fr_{X_a}v$,
\begin{equation*}
(\Fr_{X_a}E_{a,\ep}G_a(\rho)v)(\xi_a,.)=(\ep-F_\rho(\xi_a))^{-1}\hat v
\otimes \psi_{\alpha'},
\end{equation*}
so the mapping properties of $G_a(\rho)$ on $\Range E_a$ are given by
the two-body results of Weder \cite[Theorem~1.1]{Weder:Generalized},
see Remark~3.4 in \cite{Uhlmann-Vasy:Low} for its application in this
context.

For the behavior of $G_a(\rho)$ as $|z|\to\infty$, we decompose $G_a$
as above. On the range of $E_a$ we use the analogous two-body
result of Weder, while on the range of $\Id-E_a$ we use the
uniform $L^2$ operator
bounds on $R^a(\sigma)(\Id-e_a)$, valid uniformly
away from $[\Lambda_a,+\infty)$, together with $|\im F_\rho(\xi_a)|\to\infty$
for almost every $\xi_a$, namely the $\xi_a$ such that $\xi_a\cdot\nu\neq 0$,
so an application of the dominated convergence theorem yields the
desired strong convergence.
\end{proof}

Since
\begin{equation*}
P(\rho)G_a(\rho)v=(\Id+I_a G_a(\rho))v,\ v\in\calH_p,\ p>0,
\end{equation*}
we next investigate $I_a$ on $\calH_r$, again as in \cite{Uhlmann-Vasy:Low}.

\begin{lemma}\label{lemma:I_a-norm}
Suppose that $\mu>0$, and $\langle x^b\rangle^\mu
V_b\in L^\infty(X^b)$ for all $b\not\leq a$.
The multiplication operator $I_a$ is in $\bop(\calH_r,\calH_p)$ provided
that $p\leq r+\mu$, $p\leq \mu$, $r\geq -\mu$.

Moreover, there exists $C>0$
(independent of $V_b$) such that the norm of $I_a$
as such an operator is bounded by $C\max_b\sup(\langle x^b\rangle^\mu |V_b|)$.
\end{lemma}

\begin{proof}[Sketch. See \cite{Uhlmann-Vasy:Low} for details.]
Consider the matrix decomposition of $I_a$ corresponding to the direct sum in
\eqref{eq:calH_p-def}, and use the rapid decay of $\psi_{\alpha'}$ to
analyze the terms involving $E_a$. The restrictions
$p\leq r+\mu$, $p\leq \mu$, resp.\ $r\geq -\mu$ arise from the requirements of
making $E_a I_a E_a$, $E_a I_a (\Id-E_a)$, resp.\ $(\Id-E_a)I_a E_a$
bounded.
\end{proof}

Thus, perturbation theory gives

\begin{prop}\label{prop:G(rho)-exist)}
Suppose that $\rho_\perp$ satisfies $|\rho_\perp|^2+\ep_\alpha\in
(\ep_\alpha,\evth)\setminus\Lambda'_a$,
$\mu>\max(p,1)$, $p>0$, $r<0$, $r<p-1$.
There exists $\delta>0$ with the following property. Suppose that
for all $b\not\leq a$, $\sup |\langle
x^b\rangle^\mu V_b|<\delta$.
Then the operator
\begin{equation}\label{eq:G(rho)-def}
G(\rho)=G_a(\rho)(\Id+I_a G_a(\rho))^{-1}:\calH_p\to \calH_r
\end{equation}
satisfies
\begin{equation}\begin{split}\label{eq:PG-Id}
&P(\rho)G(\rho)=\Id:\calH_p\to \calH_p,\\
&G(\rho)P(\rho)=\Id:\calH_p\to \calH_p.
\end{split}\end{equation}
Moreover, $G(\rho)$ is a continuous function of $\rho$ in
$\Cx(X_a)_\alpha^\circ$, and an analytic function of $z\in\Cx\setminus\Real$,
and $\slim_{|z|\to\infty}G(\rho)=0$ as a map $\calH_p\to\calH_r$ provided that
$|z|\to\infty$ in $|\im z|\geq C|\re z|$, $C>0$.
\end{prop}

Since $I_a\psi_\alpha\in \calH_p$ for some $p>0$ if $\mu>\dim X_a/2$,
we deduce the following corollary.

\begin{cor}\label{cor:u_rho-exists}
Suppose that $\mu>\dim X_a/2$,
$|\rho_\perp|^2+\ep_\alpha\in(\ep_\alpha,\evth)\setminus\Lambda'_a$.
Then
\begin{equation}
u_\rho=u^0_\rho-e^{i\rho\cdot x_a}G(\rho)I_a \psi_\alpha=e^{i\rho\cdot x_a}
(\psi_\alpha-G(\rho)I_a \psi_\alpha)
\end{equation}
satisfies $P(\rho)u_\rho=0$, and
\begin{equation}
u^0_\rho=e^{i\rho\cdot x_a}(\Id+G_a(\rho)I_a)e^{-i\rho\cdot x_a}u_\rho.
\end{equation}
Moreover, $u_\rho-u^0_\rho\to 0$ in $\calH_r$, $r<0$,
$r<\mu-1-\frac{\dim X_a}{2}$,
as $|z|\to\infty$ in $|\im z|>C|\re z|$, $C>0$.
\end{cor}

We now only need to connect $u_\rho$ to the scattering matrix.
This is done by letting $\rho$ become real. Below we consider
\begin{equation}\label{eq:re-im-z-cond}
\im z\geq 0,\ \re z\geq 0.
\end{equation}

\begin{prop}\label{prop:G_a-limit}
The operator $G_a(\rho):\calH_p\to \calH_r$ extends continuously
(from $\Cx(X_a)_\alpha^+\cap\Cx(X_a)_\alpha^\circ$) to
$\Cx(X_a)_\alpha^+$ for $p>1/2$, $r<-1/2$, and it satisfies
$P_a(\rho)G_a(\rho)=\Id$, $G_a(\rho)P_a(\rho)=\Id$,
on $\calH_p$, $p>1/2$. The limit $G_a(\nu,\rho_\perp,z\pm i0)$ satisfies
\begin{equation}\begin{split}\label{eq:G_a-real}
&e^{i\rho\cdot x_a}G_a(\nu,\rho_\perp,z\pm i0)e^{-i\rho\cdot x_a}\\
&=R_a(\lambda+i0)+\sum_{\alpha'}\Fr^{-1}_{X_a}H(\nu\cdot(\xi_a-\re\rho))
[(|\xi_a|^2-(\lambda-\ep_{\alpha'}-i0))^{-1}\\
&\qquad\qquad\qquad\qquad\qquad
-(|\xi_a|^2-(\lambda-\ep_{\alpha'}+i0))^{-1}]\Fr_{X_a}E_{a,\alpha'},
\end{split}\end{equation}
$\lambda=\rho^2+\ep_\alpha$.
\end{prop}

\begin{proof}[Sketch of proof.]
Write
\begin{equation}\label{eq:G_a-decomp}
G_a(\nu,\rho_\perp,z)=(\Id-E_a)G_a(\nu,\rho_\perp,z)
+E_a G_a(\nu,\rho_\perp,z).
\end{equation}
The first term is just
$\Frinv_{X_a}R^a(F_\rho(.))(\Id-E_a)\Fr_{X_a}$, and $\im F_\rho\to 0$
as $\im z\to 0$. Moreover, in the region where the sign of $\im z$
matters, i.e.\ where $\re F_\rho\geq\Lambda_a$, $\im F_\rho>0$ by
\eqref{eq:sgn-im-F-p} which is valid with $\im F_\rho$ in place
of $|\im F_\rho|$ due to \eqref{eq:re-im-z-cond}, so the first term goes to
$\Frinv_{X_a}R^a(|\xi_a|^2+2\rho\cdot\xi_a+\ep_\alpha+i0)(\Id-E_a)\Fr_{X_a}$
as $\im z\to 0$. Conjugation by $e^{i\rho\cdot x_a}$ replaces $\xi_a$ by
$\xi_a+\rho$, hence yielding $R_a(\lambda+i0)(\Id-E_a)$,
$\lambda=\rho\cdot\rho+\ep_\alpha$, in the
limit.

The second term on the right hand side of \eqref{eq:G_a-decomp}
is dealt with exactly as in the two-body setting discussed by Weder.
Namely, working on the range of $E_{a,\ep}$,
$\Fr_{X_a}G_a(\rho)\Fr^{-1}_{X_a}$ is multiplication by
$(\ep-F_\rho(\xi_a))^{-1}$. For each $\xi_a$ with $\nu\cdot\xi_a\neq 0$,
the limit as $\im z\to 0$ is multiplication by the limit of
$(\ep+|\xi_a|^2+2\rho\cdot\xi_a-\ep_\alpha)^{-1}$, which is
$(\ep+|\xi_a|^2+2\rho\cdot\xi_a\pm i0-\ep_\alpha)^{-1}$, with $+$
corresponding to $\nu\cdot\xi_a>0$ and $-$ to $\nu\cdot\xi_a<0$.
Now a simple rewriting proves the proposition; see \cite{Uhlmann-Vasy:Low}
for details.
\end{proof}

Perturbation theory now shows that $G(\rho)$ itself
has a limit when $\rho$ becomes real, provided that
$\lambda=\rho^2+\ep_\alpha<\evth$. An immediate consequence is:

\begin{cor}\label{cor:u_rho-limit}
Suppose that $|\rho_\perp|^2+\ep_\alpha\in(\ep_\alpha,\evth)
\setminus\Lambda'_a$, $\mu>(\dim X_a+1)/2$.
There exists $\delta>0$ with the following property.
Suppose that
for all $b\not\leq a$, $\sup |\langle
x^b\rangle^\mu V_b|<\delta$.
Then $u_\rho$ extends
continuously to
\begin{equation}\label{eq:limit-region}
\{z:\ \im z>0\}\cup\{z\geq 0:
\ z^2+\rho_\perp^2+\ep_\alpha\in(\ep_\alpha,\evth)\setminus\Lambda'\},
\end{equation}
with $u_\rho-u^0_\rho\in \calH_r$
for all $r<-1/2$,
and $u_\rho$ is analytic in $z$ in $\im z>0$.
\end{cor}

The last part of the proof is completely analogous to the two-body
argument, expressing the S-matrix as a pairing, and considering
an analogous pairing for the complex exponential solutions.
For 2-clusters $\alpha$, $\beta$, the knowledge of $\calS_{\alpha\beta}(\lambda)$
is equivalent to that of a
renormalization $\calS^\sharp_{\alpha\beta}(\lambda)$
with free scattering removed, see \cite[Section~2]{Uhlmann-Vasy:Low}.
Rather than giving its detailed definition here, we simply give an expression
for its Schwartz kernel that is derived in
\cite[Corollary~2.2]{Uhlmann-Vasy:Low}.
For
$\rho\in X_a$ (real!), we write
\begin{equation}\begin{split}
&U_\rho=(\Id-R(\lambda+i0)I_a)u^0_\rho,\ \lambda=\rho^2+\ep_\alpha,\\
&u^0_\rho=u^0_{\alpha,\rho}=e^{i\rho\cdot x_a}\psi_\alpha(x^a).
\end{split}\end{equation}
Then, for $\lambda\nin\Lambda'$, $\langle x^b\rangle^\mu V_b\in L^\infty(X^b)$,
$\mu>\dim X_a$,
\begin{equation}\label{eq:pair-8}
\calS^\sharp_{\alpha\beta}
(\lambda,\omega,\omega')=\int_{\Rn} I_b U_\omega(x)
\overline{e^{ix_b\cdot\omega'}\psi_\beta(x^b)}\,dx=\int_{\Rn} I_b U_\omega(x)
e^{-ix_b\cdot\omega'}\overline{\psi_\beta(x^b)}\,dx;
\end{equation}
in particular, $\calS^\sharp_{\alpha\beta}(\lambda)$ has a continuous kernel.
The analogous pairing for the complex exponential solutions is
\begin{equation}\begin{split}\label{eq:pair-16}
G_{\alpha\alpha'}(\rho,\overline{\rho}+\zeta)
&=\int_{\Rn} I_a(x) u_{\rho}(x)\overline{u^0_{\alpha',\overline{\rho}+\zeta}(x)}\,dx\\
&=\int_{\Rn} I_a(x) u_\rho(x) \overline{\psi_{\alpha'}(x^a)}
e^{-i(\rho+\zeta)\cdot x_a}\,dx.
\end{split}\end{equation}
By Corollary~\ref{cor:u_rho-limit}, if
$\langle x^b\rangle^\mu V_b\in L^\infty(X^b)$ for all $b$
and for some $\mu>\dim X_a$, then the integral in \eqref{eq:pair-16}
converges for all $\rho$ for which $u_\rho$ exists, and for
$\zeta\in\Rn$, since then the real parts of the
exponentials cancel, and
$I_a \psi_\alpha\overline{\psi_{\alpha'}}\in L^1(X_0)$,
and the same holds for $I_a \overline{\psi_{\alpha'}} G(\rho)(I_a\psi_\alpha)$.
Other properties of \eqref{eq:pair-16} follow immediately from
Corollary~\ref{cor:u_rho-limit}.

\begin{prop}
Suppose that $|\rho_\perp|^2+\ep_\alpha\in(\ep_\alpha,\evth)
\setminus\Lambda'_a$, $\mu>\dim X_a$,
and $V_b$ as in
Corollary~\ref{cor:u_rho-limit}. Then
$G_{\alpha\alpha'}$ is an
analytic function of $z$ in $\Cx\setminus\Real$, and extends to be
continuous on \eqref{eq:limit-region}.
In addition,
\begin{equation}\label{eq:G-z-infty}
\lim_{|z|\to\infty} G_{\alpha\alpha'}(\rho,\overline{\rho}+\zeta)
=\int_{\Rn} I_a \psi_\alpha \overline{\psi_{\alpha'}} e^{-i\zeta\cdot x_a}\,dx,
\end{equation}
provided that $|z|\to\infty$ in $|\im z|>C|\re z|$, $C>0$.
\end{prop}

For fixed $\rho$ real and $\zeta$ satisfying
\begin{equation*}
\lambda=\rho^2+\ep_\alpha=(\rho+\zeta)^2+\ep_{\alpha''},
\end{equation*}
i.e.\ the equality of incoming and outgoing energies, we can
relate $G_{\alpha\alpha''}(\rho,\rho+\zeta)$, $\rho$ real, to the
S-matrices as follows.
Under our assumptions,
\begin{equation}\label{eq:R-G-8}
e^{i\rho\cdot x_a}(\Id+G_a(\rho)I_a)e^{-i\rho\cdot x_a}u_\rho
=u^0_\rho=(\Id+R_a(\lambda+i0)I_a) U_\rho.
\end{equation}
Applying $(\Id+R_a(\lambda+i0)I_a)^{-1}$ to both sides of
\eqref{eq:R-G-8}, and using \eqref{eq:G_a-real}, we deduce that
\begin{equation}\begin{split}\label{eq:R-G-16}
U_\rho&=u_\rho\\
&+2\pi i\sum_{\alpha'}
(\Id+R_a(\lambda+i0)I_a)^{-1}\Fr^{-1}_{X_a}H(\nu\cdot(\xi_a-\re\rho))
\delta_{|\xi_a|^2-(\lambda-\ep)}
\Fr_{X_a}E_{a,\alpha'}I_a u_\rho.
\end{split}\end{equation}
Here we used
\begin{equation*}
(\rho-(\lambda-\ep-i0))^{-1}
-(\rho-(\lambda-\ep+i0))^{-1}
=-2\pi i \delta_{\rho-(\lambda-\ep)}
\end{equation*}
with $\rho=|\xi_a|^2$; the $\delta$ distribution here is the
pull-back of the standard delta distribution on $\Real$ by the
map $\xi_a\mapsto |\xi_a|^2-(\lambda-\ep)$.
Integrating against
$I_a e^{-i(\rho+\zeta)\cdot x_a}\overline{\psi_{\alpha''}(x^a)}$, and
yields
\begin{equation}\begin{split}
\calS^\sharp_{\alpha\alpha''}(\lambda,\rho,\rho+\zeta)
=&G_{\alpha\alpha''}(\rho,\rho+\zeta)\\
&-\sum_{\alpha'}
\frac{i}{2\sqrt{\lambda-\ep_{\alpha'}}}
\int_{\sphere_a(\sqrt{\lambda-\ep_{\alpha'}})}
\calS^\sharp_{\alpha'\alpha''}(\lambda,\rho',\rho+\zeta)\\
&\qquad\qquad\qquad\qquad H(\nu\cdot (\rho'-\re\rho))
G_{\alpha\alpha'}(\rho,\rho')\,d\rho'.
\end{split}\end{equation}
This is an integral equation for $G_{\alpha\alpha''}$ in terms of
$\calS^\sharp_{\alpha\alpha''}(\lambda)$,
$\calS^\sharp_{\alpha'\alpha''}(\lambda)$,
$\lambda=\rho^2+\ep_\alpha$. It can be solved due to the smallness
assumption on $I_a$, so we deduce the following result, with the
same proof as in \cite{Uhlmann-Vasy:Low}.

\begin{prop}\label{prop:G-S}(\cite[Proposition~5.2]{Uhlmann-Vasy:Low})
Suppose that $\rho^2+\ep_\alpha=\lambda
=(\rho+\zeta)^2+\ep_{\alpha''}\in(-\infty,0)\setminus\Lambda'$,
$|\rho_\perp|^2+\ep_\alpha\in(\ep_\alpha,\evth)
\setminus\Lambda'_a$, $\mu>\dim X_a$. There exists
$\delta>0$ with the following property.

Suppose that for all $b\not\leq a$, $\sup |\langle x^b\rangle^\mu V_b|<\delta$.
Then the pairings $G_{\alpha\alpha''}(\rho,\rho+\zeta)$
are determined by the operators
$\calS^\sharp_{\alpha'\alpha''}(\lambda)$ given for all
$\alpha'$ and $\alpha''$.
\end{prop}

We can now prove Theorem~\ref{thm:main}.
Fix a non-empty open interval $I\subset(\ep_\alpha,\inf\Lambda_a)$
with $\sup I<\inf\Lambda_a$, and let
$R=2\sqrt{\sup I-\ep_\alpha}$. Let $\zeta\in X_a$ satisfy $|\zeta|<R$, and
$\frac{1}{4}|\zeta|^2+\ep_\alpha\nin\Lambda_a'$. The last condition
excludes a discrete set of values of $|\zeta|^2$ in $[0,R)$.
Note that $\frac{1}{4}|\zeta|^2+\ep_\alpha<0$.

We now choose $\evpr<\inf\Lambda_a$ in \eqref{eq:ep_1-def} so that
$\frac{1}{4}|\zeta|^2+\ep_\alpha<\evpr$, and
\begin{equation}\label{eq:I-energy}
I\cap(\frac{1}{4}|\zeta|^2+\ep_\alpha,\evpr)\neq\emptyset.
\end{equation}
Let $\rho_\perp=-\zeta/2$, so $\rho_\perp^2+\ep_\alpha\in(\ep_\alpha,\evpr)
\setminus\Lambda'_a$.
Let $\nu\in X_a$ be such that $|\nu|^2=1$, $\nu\cdot \zeta=0$.
Then for $\rho=z\nu+\rho_\perp$, $z\in\Cx$,
$(\rho+\zeta)^2-\rho^2=2\rho\cdot\zeta+\zeta^2=0$.
Having fixed $\zeta$, $\rho_\perp$, $\nu$, consider the
energy $\rho^2+\ep_\alpha$ when $z$ becomes real.
Since
$\rho^2=z^2+|\rho_\perp|^2$, as $z$ varies over
$[0,\sqrt{\evth-\rho_\perp^2-\ep_\alpha})$, the
energy varies over $[|\rho_\perp|^2+\ep_\alpha,\evth)$. This
intersects the given interval $I$ by \eqref{eq:I-energy} as $\evpr<\evth$.

Let $z_0>0$ satisfy $\lambda=\rho_\perp^2+z_0^2+\ep_\alpha
\in(\ep_\alpha,\evth)\setminus\Lambda'$.
By Proposition~\ref{prop:G-S}
the limit of the pairing $G_{\alpha\alpha}(\rho,\rho+\zeta)$
as $z\to z_0$ (in $\im z>0$)
is determined by the scattering matrices
$\calS^\sharp_{\alpha'\alpha''}(\rho_\perp^2+z_0^2+\ep_\alpha)$. Thus,
there exists a non-empty open interval of these values $z_0$ at which
$G_{\alpha\alpha}(\rho,\rho+\zeta)$ is determined
by $\calS^\sharp_{\alpha'\alpha''}(\lambda)$, $\lambda\in I$.

Since the limit on any open interval in the boundary of its domain
determines an analytic function, we deduce that knowing the S-matrix
$\calS^\sharp_{\alpha'\alpha''}$ in the interval $I$
determines $G_{\alpha\alpha}(\rho,\zeta)$
for all $\zeta$ with $|\zeta|<R$.

Now let $z\to\infty$ through imaginary $z$.
By \eqref{eq:G-z-infty}, $G_{\alpha\alpha}(\rho,\rho+\zeta)$ converges to
\begin{equation}\label{eq:pair-99}
\int_{\Rn} I_a |\psi_{\alpha}|^2 e^{-i\zeta\cdot x_a}
\,dx=\int_{X_a}e^{-i\zeta\cdot x_a} \left(\int_{X^a}
I_a|\psi_{\alpha}|^2\,dx^a\right)\,dx_a.
\end{equation}
But this is the Fourier transform of $\int_{X^a}
I_a|\psi_{\alpha}|^2\,dx^a$ in $X_a$, evaluated at $\zeta$.
Hence $\calS_{\alpha'\alpha''}(\lambda)$, $\lambda\in I$, determines
the Fourier transform in $x_a$ of the `effective interaction'
\begin{equation}\label{eq:pair-101}
\int_{X^a}
I_a|\psi_{\alpha}|^2\,dx^a
\end{equation}
in a ball of radius $2\sqrt{\sup I-\ep_\alpha}$, except possibly on the spheres
$\frac{1}{4}|\zeta|^2+\ep_\alpha\in\Lambda'_a$. However, $\int_{X^a}
I_a|\psi_{\alpha}|^2\,dx^a\in L^1(X_a)$, hence the Fourier transform is
continuous, hence it is determined on the whole ball.
This completes the proof of Theorems~\ref{thm:main} and
\ref{thm:main2}.

Finally, we return to the electron-ion scattering experiment, comparing
its treatment as a two-body problem and as a many-body problem.
If one only uses a fixed energy $\lambda$, and the potentials are
exponentially decaying, in the two-body model the S-matrix determines
the interaction, which is $V_{\alpha,\eff}$ in this model
(here $\alpha$ is the ground state
of the ion). No such result exists if we consider the ion as a composite
particle, i.e.\ if the experiment is treated as a many-body problem,
although our theorem shows that if the 2-cluster to 2-cluster S-matrix is
known in an {\em interval} of energies and if the interactions are weak,
then we can determine
$V_{\alpha,\eff}$. Indeed, in the two-body fixed energy
inverse result one has to let $\rho_\perp\to\infty$. If we consider
the scattering experiment as a
many-body problem, $\rho_\perp\to\infty$
causes serious complications, namely the possibility
of the break up of the ion. This seems to limit the use of two-body
results for composite particles, such as ions.

\bibliographystyle{plain}
\bibliography{sm}

\end{document}